\documentclass[times,sort&compress,3p]{elsarticle}
\journal{Journal of Multivariate Analysis}
\usepackage[labelfont=bf]{caption}

\usepackage{amsmath,amsfonts,amssymb,amsthm,booktabs,color,epsfig,graphicx,hyperref,url,enumerate}
\usepackage{dsfont}

\newcommand{\1}{\mathds{1}}
\newcommand{\diag}{\mathrm{diag}}
\newcommand{\etr}{\mathrm{etr}}
\newcommand{\sgn}{\mathrm{sgn}}
\newcommand{\tr}{\mathrm{tr}}

\newcommand{\eg}{{\it e.g.}}
\newcommand{\etal}{\textit{et al.}}
\newcommand{\ie}{{\it i.e.}}

\newcommand{\pd}[1]{\partial_{#1}}
\newcommand{\qq}{q}

\newtheorem{theorem}{Theorem}
\newtheorem{proposition}{Proposition}
\newtheorem{lemma}{Lemma}
\newtheorem{example}{Example}
\newtheorem{heuristic}{Heuristic}

\begin{document}

\begin{frontmatter}

\title{Computation of the expected Euler characteristic for the largest eigenvalue of a real non-central Wishart matrix}

\author[A1]{Nobuki Takayama}
\author[A2]{Lin Jiu}
\author[A3]{Satoshi Kuriki}
\author[A4,A5]{Yi Zhang \corref{mycorrespondingauthor}}

\address[A1]{Department of Mathematics, Kobe University, Japan}
\address[A2]{Department of Mathematics and Statistics, Dalhousie University, Canada}
\address[A3]{The Institute of Statistical Mathematics, Research Organization of Information and Systems, Japan}
\address[A4]{Department of Mathematical Sciences, The University of Texas at Dallas, USA}
\address[A5]{Department of Mathematical Sciences, Xi'an Jiaotong-Liverpool University, China}

\cortext[mycorrespondingauthor]{Corresponding author. Email address: \url{Yi.Zhang03@xjtlu.edu.cn}}

\begin{abstract}
We give an approximate formula for the distribution of the largest
eigenvalue of real Wishart matrices by the expected Euler characteristic
method for general dimension.
The formula is expressed in terms of a definite integral with parameters.
We derive a differential equation satisfied by the integral
for the $2\times 2$ matrix case and perform a numerical analysis of it.
\end{abstract}

\begin{keyword} 
Euler characteristic method\sep
holonomic gradient method\sep
real non-central Wishart distributions.
\MSC[2010] 62H10 \sep
68W30
\end{keyword}
\end{frontmatter}

\section{Introduction}
\label{sec:INTRO}

For $i \in \{1,\ldots,n\}$, let $\xi_i\in\mathbb{R}^{m\times 1}$ be independently distributed as the $m$-dimensional (real) 
Gaussian distribution $\mathcal{N}_m(\mu_i,\Sigma)$, where $\mu_i$ and $\Sigma$ 
are the mean vector and covariance matrix of $\xi_i$, respectively.
The (real) Wishart distribution $\mathcal{W}_m(n,\Sigma;\Omega)$ is the probability measure on the cone of $m\times m$ 
positive semi-definite matrices induced by the random matrix
\[
 W = \Xi \Xi^{\top}, \quad \Xi=(\xi_1,\ldots,\xi_n)\in\mathbb{R}^{m\times n}.
\]
Here, $\Omega=\Sigma^{-1}\sum_{i=1}^n \mu_i \mu_i^{\top}$ is a parameter matrix.
Unless $\Omega$ vanishes, the corresponding distribution is referred to as the non-central (real) Wishart distribution.

The largest eigenvalue $\lambda_1(W)$ of $W$ is used as a test statistic for testing $\Sigma=I_m$ and/or $\Omega\ne 0$ 
under the assumption that $\Sigma-I_m$ is positive semi-definite.
This test statistic is expected to have good power when the matrices $\Sigma-I_m$ and $\Omega$ are of  small size.

When testing hypotheses, the distribution of $\lambda_1(W)$, which is the largest eigenvalue of $W$, is of particular interest as it gives the power of the test.
When $\Omega=0$, the works by James and other authors 
(see, \eg, Muirhead~\cite{Muirhead2005}) show
that the cumulative distribution function of $\lambda_1(W)$ can be written as a hypergeometric function with matrix argument as follows:
\begin{equation*}
 \Pr(\lambda_1(W)<x) = c_{m,n} \det\left(\frac{1}{2} n x \Sigma^{-1}\right)^{n/2}{}_1 F_1\left(\frac{1}{n};\frac{1}{2}(n+m+1);-\frac{1}{2}n x\Sigma^{-1}\right),
\end{equation*}
where $c_{m,n}$ is a known constant \cite[Corollary 9.7.2]{Muirhead2005}.
It is well known that the hypergeometric function ${}_1 F_1$ has a series expression in the zonal polynomial $C_\kappa$ 
with index $\kappa$, which is a partition of an integer.
However, in view of numerical calculation, this is less useful because the explicit form of $C_\kappa(X)$ is not known 
unless the rank of the matrix $X$ is 1 or 2.
On account of this difficulty, Hashiguchi \etal~\cite{Hashiguchi2013} proposed a holonomic gradient method (HGM) for numerical evaluation,
which utilizes a holonomic system of differential equations 
for computation.
However, when $\Omega\ne 0$, the situation is more difficult. In this case, the cumulative distribution function $\Pr(\lambda_1(W)<x)$ cannot be expressed as a simple series of  zonal polynomials.
Hayakawa~\cite[Corollary 10]{Hayakawa1969} provides a formula for the cumulative distribution function as a series expansion in the Hermite polynomial $H_\kappa$ with symmetric matrix argument defined by the Laplace transform of~$C_\kappa$ as follows:
\[
 \etr\bigl(-T T^{\top}\bigr) H_\kappa(T) = \frac{(-1)^{|\kappa|}}{\pi^{m n/2}}
 \int \etr\bigl(-2i T U^{\top}\bigr) \etr\bigl(-U U^{\top}\bigr) C_\kappa\bigl(U U^{\top}\bigr)\,dU,
 \quad T, U\in\mathbb{R}^{m\times n}.
\]
The Hermite polynomial $H_\kappa$ can be written as a linear combination of the zonal polynomial $C_\kappa$;  however,  the coefficients 
 not provided explicitly \cite{Chikuse1992}.
Another approach is  to use invariant polynomials proposed by Davis~\cite{Davis1979,Davis1980}.
Using the probability density function of the non-central Wishart distribution derived by James \cite{James1955},
the cumulative distribution function of $\lambda_1(W)$ is shown to be proportional to
\[
 \int_{0<W<xI_m} 
 |W|^{(n-m+1)/2-1}\etr\left(-\frac{1}{2}\bigl(\Sigma^{-1}W+\Omega\bigr)\right) {}_0F_1(n/2;\Omega\Sigma^{-1}W/4)\,dW.
\]
D\'{i}az-Gar\'{i} and Guti\'{e}rrez-J\'{a}imez 
 \cite{Garcia2011}  showed that this has a series expansion in terms of invariant polynomials.
Here, the invariant polynomial is a polynomial in two matrices indexed by two partitions.
Although, in principle, the invariant polynomial can be expressed in terms of zonal polynomials in two matrices, it is  challenging to utilize this expression for numerical calculation.

In this paper, instead of the direct calculation approach, we approximate the distribution function through the 
expected Euler characteristic heuristic or the Euler characteristic method proposed  by Adler \cite{adler-1981} and Worsley \cite{worsley-1995}.
(see also \cite{Adler2007} and \cite{kuriki-takemura-2009}.)
This is a methodology to approximate the tail upper probability of a random field.
In our problem, since the square root of the largest eigenvalue $\lambda_1(W)^{1/2}$ is the maximum of a Gaussian field
\[
 \bigl\{ u^{\top} \Xi v \mid \Vert u\Vert_{\mathbb{R}^m}=\Vert v\Vert_{\mathbb{R}^n}=1 \bigr\},
\]
this method  is applicable (~\cite{kuriki-takemura-2001},~\cite{kuriki-takemura-2008}).
One can show that the Euler characteristic method evaluates the quantity
\begin{equation}
\label{altsum}
 \Pr(\lambda_1(W)\ge x) - \Pr(\lambda_2(W)\ge x) + \cdots + (-1)^{m-1} \Pr(\lambda_m(W)\ge x)
\end{equation}
rather than $\Pr(\lambda_1(W)\ge x)$.
Nevertheless, this formula approximates $\Pr(\lambda_1(W)\ge x)$ well when $x$ is large because $\Pr(\lambda_i(W)\ge x)$ ($i\ge 2$) 
are negligible for large $x$.
This is practically sufficient for our purpose because only the upper tail probability is required in testing hypotheses.

In this paper, we consider the non-central real Wishart matrix.
In the multiple-input multiple-output (MIMO) problem, the non-central complex Wishart matrix also plays an important role.
The largest eigenvalue of the non-central complex Wishart  matrix is  significantly easier to handle in that case
because the explicit formula for the cumulative distribution  was given
by Kang and Alouini \cite{Kang2003}.
The holonomic gradient method based on Kang and Alouini's formula was proposed 
in~\cite{Danufane2017}.

In general, the approximation error of the Euler characteristic has not been extensively studied.
Nevertheless, the Euler characteristic heuristic is widely used as an approximation of the tail probability of the supremum because of the difficulty of the original problem and the empirical usefulness of this heuristic (see, \eg,~\cite{taylor-worsley-2013}).
One exception is the Gaussian process with mean zero and variance one,
which corresponds to the central Wishart case where $\Sigma$ is proportional to the identity matrix and $\Omega$ vanishes.
In this particular case, the approximation error has been fully investigated~\cite{kuriki-takemura-2001,kuriki-takemura-2008}; 
the details are presented in~\ref{sec:central}.
For the non-central case, we present the following lemma for which the proof is provided in~\ref{SEC:lemma2}:
\begin{lemma}
\label{lem:m=2} 
Assume that $m=2$. If either $\Sigma$ or $\Sigma\Omega$ has distinct eigenvalues, then it holds that
\[
 \Pr(\lambda_2(W)\ge x) = o\Bigl(\Pr(\lambda_1(W)\ge x\Bigr) \quad \mbox{as $x\to\infty$}.
\]
\end{lemma}
This implies that the Euler characteristic approximation (\ref{altsum}) is justified as  an approximation for $\Pr(\lambda_1(W)\ge x)$.
We conjecture that this holds for arbitrary $m$ and arbitrary configurations of $\Sigma$ and $\Omega$.

The rest of the paper is organized as follows.
In Section \ref{sec:expectation}, we provide  an integral representation
formula for the expectation of the Euler characteristic for random matrices
of a general size.
In Section \ref{sec:case}, we consider the case of $2 \times 2$ random
matrices and study the integral representation derived 
in Section \ref{sec:expectation} in the polar coordinate system
and investigate it from a numerical point of view.
By virtue of the theory of holonomic systems
(see, \eg, \cite{Hibi2013}), the integral representation
given in Section \ref{sec:expectation} satisfies a holonomic system
of linear differential equations. 
However, its explicit form is not known in general.
In Section \ref{sec:computer}, we consider the case of $2\times 2$ random
matrices again to demonstrate that
the recent development~\cite{KKZ2009, KKZ2011} of computer aided proofs and derivations 
(CAPD),
for combinatorial identities, proofs of them, derivations of 
difference, and differential equations
can be applied to the evaluation of definite integrals or sums.
We derive a differential equation which is satisfied by the integral representation
of the expectation of the Euler characteristic
with the help of computer algebra algorithms
and perform a numerical analysis of the differential equation.
Thus, a new efficient method to numerically evaluate the Euler expectation,
when the numerical integration is difficult to perform, is obtained. Last but not least, 
in~\ref{sec:central},
 we give a closed formula, expressed in terms of the Laguerre polynomial, for the expectation of the Euler characteristic for random matrices
of general size for the central and scalar covariance case.

\section{Expectation of an Euler characteristic number}
\label{sec:expectation}

Let $A=(a_{ij})$ be a real $m \times n$ matrix-valued random variable (random matrix)
with density 
$$ p(A) dA, \quad dA=\prod da_{ij}. $$
We assume that $p(A)$ is smooth and $n \geq m \geq 2$.
Define a manifold 
$$ M=\{ h g^{\top} \,|\, g \in S^{m-1}, h \in S^{n-1} \} \simeq S^{m-1}\times S^{n-1}/\sim, $$
where $(h,g) \sim (-h,-g)$, $h$ and $g$ are column vectors,
and $h g^{\top}$ is a rank $1$ $m \times n$ matrix.
Set
$$ f(U) = \tr(UA)= g^{\top} A h, \quad U \in M, $$
and 
$$ M_x = \{ h g^{\top} \in M \,|\, f(U)=g^{\top} A h \geq x \}. $$

\begin{proposition}
Let $A$ be a random matrix as  aforementioned.
The following claims are equivalent:
\begin{enumerate}
\item[(i)] The function $f(U)$ has a critical point at $U=hg^{\top}$.
\item[(ii)] The vectors $g^{\top}, h$ are left and right eigenvectors of $A$, respectively.
In other words, there exists a constant $c$ such that
$g^{\top} A = c h^{\top}$, $Ah = cg$. 
\end{enumerate}
Moreover, the function $f$ takes value $c$ at the critical point $(g,h)$.
\end{proposition}

\textbf{Proof.}
We assume that $g \in S^{n-1}$ and $h \in S^{m-1}$ are expressed by local coordinates $u_i$ and $v_a$, respectively, 
where $1\leq i \leq m-1$ and $1 \leq a \leq n-1$.
We denote $\partial/\partial u_i$ by $\pd{i}$
and $\partial/\partial v_a$ by $\pd{a}$.
Since $g^{\top} g = 1$, we have
$g_i^{\top} g = 0$, where $g_i = \pd{i}\bullet g$.
We omit $\bullet$, which  represents the action, when there is no ambiguity.
Analogously, we have $h_a^{\top} h_a = 0$, where $h_a = \pd{a} h$.

Assume that $A$ is a $m \times n$ (random real) matrix.
Let us consider the function $f(U)$ expressed by the local coordinate $(g(u),h(v))$
\[
 f(g,h) = g^{\top} A h, \quad g \in S^{n-1}, \ \ h \in S^{m-1}.
\]
At the critical point of $f$, we have
\begin{equation}
\label{EQ:equality}
\pd{i} f = g_i A h = 0, \quad
  \pd{a} f = g A h_a = 0.
\end{equation}
The aforementioned equality~\eqref{EQ:equality}holds for each $i$, and $u$ is a local coordinate of $S^{n-1}$, which implies that
all $g_i$'s are linearly independent.
Therefore, there exists a constant $c$ such that $Ah = cg$ at the critical point.
Analogously, there exists a constant $d$ such that
$A^{\top} g = d h$.
Let us show that $c=d$.
We have
$$ (g^{\top} A) h = (d h^{\top}) h = d (h^{\top} h) = d $$
and
$$ g^{\top} (A h) = g^{\top} (cg) = c (g^{\top} g) = c. $$
Therefore, we have $d=c=f(g, h)$ at the critical point.

Conversely, $Ah=cg$ and $A^{\top} g = dh$ at a point $(u,v)$ imply that
$(g(u),h(v))$ is a critical point of $f(g(u),h(v))$.
\qedsymbol

We consider a continuous family of elements of  $\mathit{SO}(m)$ parameterized by the first column
vector $g$.
In other words, we consider a continuous family of orthogonal frames of $\mathbb{R}^m$
parameterized by $g \in S^{m-1}$.
 An element of $\mathit{SO}(m)$ is denoted by $(g, G) \in O(m)$, 
where $G$ is an $m \times (m-1)$ matrix.
Analogously, we take a family $(h,H) \in \mathit{SO}(n)$
parameterized by $h \in S^{n-1}$ , 
where $H$ is an $n \times (n-1)$ matrix parameterized by $h$.
Set
\[
  \sigma = g^{\top} A h, \ 
  B = G^{\top}(g) A H(h).
\]
The matrix $A$ can be expressed as
\begin{equation}
\label{eq:ghB:coordinate}
  A = \sigma g h^{\top} + G(g) B H(h)^{\top}.
\end{equation}
Intuitively, this is a partial singular value decomposition.
We denote the set of  all $(m-1) \times (n-1)$ matrices by
$M(m-1,n-1)$.

 The aforementioned decomposition provides a coordinate system for the space of random matrices~$A$.
Without loss of generality, we assume that $m \leq n$. 
We sort the singular values of $B$  in descending order, and denote by $\lambda_j(B)$ the $j$-th singular value of the matrix $B$.
For a real number $\sigma$, we define
\begin{align*}
  {\cal B}(i,\sigma) = \{ B \in M(m-1,n-1)\,|\, 
 & \mbox{all the singular values of $B$ are different and non-zero}, \\
 & \mbox{$\lambda_j(B) > \sigma$ for all $j<i$}, \ 
   \mbox{$\lambda_j(B) \leq \sigma$ for all $j\geq i$} \}. 
\end{align*}
 Subsequently, we set
$$
{\cal A} = \{ A \in M(m,n) \,|\, \mbox{all the singular values of $A$ are different and non-zero} \},
$$
and
$$
{\cal A}_i = \{ (\sigma, g, h, B) \,|\,
  \sigma \in \mathbb{R}_{>0},
 (g,h) \in S^{m-1} \times S^{n-1}/\sim,
 B \in {\cal B}(i,\sigma) \}.
$$
For a matrix $A$ in ${\cal A} \subset M(m,n)$,
we sort the singular values of $A$  in descending order as follows:
$$ \sigma^{(1)} > \sigma^{(2)} > \cdots > \sigma^{(m)} > 0. $$
Let $g^{(i)}$ and $h^{(i)}$ be the left and right eigenvectors of $A$ for $\sigma^{(i)}$, respectively.
Note that $g^{(i)}$ and~$h^{(i)}$ are respective eigenvectors of $A A^{\top}$ and $A^{\top} A$ for the eigenvalue $\sigma^{(i)}$, which 
implies that $g^{(i)}$ and~$h^{(i)}$ are uniquely determined modulo the multiplication by
$\pm 1$.  
Define a map $\varphi_i$ from ${\cal A}$ to ${\cal A}_i$ by
\begin{equation}
\label{EQ:phi}
 \varphi_i(A) = (\sigma^{(i)}, g^{(i)}, h^{(i)}, G(g^{(i)}) A H^{\top}(h^{(i)})).
\end{equation}
The matrix $G(g^{(i)}) A H^{\top}(h^{(i)})$ lies in
${\cal B}(i,\sigma^{(i)})$ because the singular values of $B^{(i)}$ agree with those of $A$ excluding
$\sigma^{(i)}$.

\begin{lemma}
\label{lemma:one-to-one-2017-12-08}
The map $\varphi_i$ in~\eqref{EQ:phi} is smooth and isomorphic.
\end{lemma}

\textbf{Proof.}
Define a map $\psi$ from ${\cal A}_i$ to ${\cal A}$ by
$$
 \psi(\sigma,g,h,B) = g \sigma h^{\top} + G(g) B H(h)^{\top}. 
$$
Based on calculation, we observe that $\varphi_i \circ \psi$ and $\psi \circ \varphi_i$
are identity maps.
The map $\varphi_i$ is then one-to-one and  surjective.
Next, we show that the map $\psi$ is smooth.
Since we assume that all the singular values are different,
the maps of taking  the $i$-th singular value of a given $A$ and an eigenvector for the singular value
are smooth on an open connected neighborhood $W \subset {\cal A}$ of $A$ (by checking the Jacobian does not vanish).
The inverse map is then locally smooth.
Hence, $\varphi_i$ is smooth and  isomorphic.
\qedsymbol

We are interested in the Euler characteristic of $M_x$.

\begin{theorem}
\label{th:expectation_of_euler_characteristic}
Suppose that $x > 0$ and $f(U)$ is a Morse function
for almost all $A$'s.
We further assume that 
if a set is  of measure zero with respect to the Lebesgue measure, it is also a measure zero set with respect to the measure $p(A) dA$.
The expectation of the Euler characteristic number ${\rm E}[\chi(M_x)]$  equals
\begin{equation}
\frac{1}{2} \int_x^\infty \sigma^{n-m} d\sigma \int_{\mathbb{R}^{(m-1)(n-1)}} dB \int_{S^{m-1}} G^{\top} dg
 \int_{S^{n-1}} H^{\top} dh \ 
 \det\left(\sigma^2 I_{m-1}-B B^{\top}\right) p(A).
\label{EchiMx}
\end{equation}
Here, we set
$G^{\top} dg = \wedge_{i=1}^{m-1} G_i^{\top} dg$,
$H^{\top} dh = \wedge_{i=1}^{n-1} H_i^{\top} dh$,
where $G_i$ and $H_i$ are the $i$-th column vectors of $G$ and $H$, respectively, $dg = (dg_1, \ldots, dg_m)^{\top}$ and
$dh = (dh_1, \ldots, dh_n)^{\top}.$
\end{theorem}
Note that $G^{\top} dg$ and $H^{\top} dh$ are $O(m)$ and $O(n)$ invariant measures on $S^{m-1}$
and $S^{n-1}$, respectively.

\textbf{Proof.}
Without loss of generality, we assume that $m \leq n$.
According to Morse theory, if $f(U)$ is a Morse function,
which is a smooth function without a degenerated critical point,
then we have
\begin{align}
 \chi(M_x) =& \sum_{\mbox{critical point}} \1(f(U) \geq x) \,\sgn\,\det
\left(
\begin{array}{cc}
 -\pd{i}\pd{j} f & -\pd{i}\pd{a} f \\
 -\pd{a}\pd{i} f & -\pd{a}\pd{b} f \\
\end{array}
\right)
\label{eq:det1} \\
=& \sum_{\mbox{eigenvectors}} \1(\sigma \geq x) \,\sgn\,\det
\left(
\begin{array}{cc}
 \sigma I_m & -GBH^{\top} \\
 -H B^{\top} G^{\top} & \sigma I_n \\
\end{array}
\right)
\label{eq:det2} \\
=& \sum_{i=1}^m \1(\sigma^{(i)} \geq x) \,\sgn\,
{\sigma^{(i)}}^{n-m} {\sigma^{(i)}}^2
\det
\left( 
{\sigma^{(i)}}^2 I_{m-1} - B^{(i)} {B^{(i)}}^{\top}
\right),
\label{eq:det3} 
\end{align}
where $\sigma^{(i)}$ is the $i$-th singular value of $A$,
$g^{(i)}$ and $h^{(i)}$ are left and right eigenvectors, and 
$B^{(i)} = G^{\top}(g^{(i)}) A H(h^{(i)})$.
The equality (\ref{eq:det1}) is the Morse theorem for manifolds with boundaries.
The equalities (\ref{eq:det1}) and (\ref{eq:det2}) can be established as follows.

First, we have the relation $g_i^{\top} g = 0$.
By differentiating it with respect to $u_j$, we have
$g_{ij}^{\top} g + g_i^{\top} g_j = 0$.
Let us evaluate $\pd{i}\pd{j} f$.
By the expression
$A=\sigma g h^{\top} + GBH^{\top}$, it is equal to 
\[
\pd{i}\pd{j} f = g_{ij}^{\top} A h = g_{ij}^{\top} \sigma g h^{\top} h + g_{ij}^{\top} GBH^{\top} h = -\sigma g_i^{\top} g_j \ \ \text{ since } \ \ H^{\top} h=0.
\]
%
Next, we evaluate $\pd{i}\pd{a} f$:
\[
\pd{i}\pd{a} f = g_{i}^{\top} A h_a = g_{i}^{\top} g \sigma h^{\top} h_a + g_{i}^{\top} GBH^{\top} h_a = g_{i}^{\top} GBH^{\top} h_a \ \ \text{ since } \ \ g_i^{\top} g = h^{\top} h_a=0.
\]
Third, we evaluate $\pd{a}\pd{b} f$:
\[
\pd{a}\pd{b} f = g^{\top} A h_{ab} = g^{\top} g \sigma h^{\top} h_{ab} + g^{\top} GBH^{\top} h_{ab} = -\sigma h_a^{\top} h_b \ \ \text{ since } \ \ g^{\top} G=0.
\]
%
Summarizing the aforedescribed calculations, the Hessian is equal to
\begin{align*}
& \left(
\begin{array}{cc}
-\pd{i}\pd{j} f & -\pd{i}\pd{a} f \\
-\pd{i}\pd{a} f & - \pd{a}\pd{b} \\
\end{array} \right) 
= \left(
\begin{array}{cc}
\sigma g_i^{\top} g_j & -g_i^{\top} GBH^{\top} h_a \\
-h_a^{\top} H B^{\top} G^{\top} g_i & \sigma h_a^{\top} h_b \\
\end{array} \right) \\
&= 
\left(
\begin{array}{cc}
g_1 \times \cdots \times g_{n-1} & 0 \\
0 & h_1 \times \cdots \times h_{m-1} \\
\end{array} \right)^{\top}
\left(
\begin{array}{cc}
 \sigma I_{m} & -GBH^{\top} \\
- H B^{\top} G^{\top} & \sigma I_{n} \\
\end{array} \right)
\left(
\begin{array}{cc}
g_1 \times \cdots \times g_{n-1} & 0 \\
0 & h_1 \times \cdots \times h_{m-1} \\
\end{array} \right).
\end{align*}
Since $\det(P P^{\top}) = \det(P)^2$,
the sign of the determinant of the Hessian is equal to that of the middle of the above 3 matrices.

The equalities (\ref{eq:det2}) and (\ref{eq:det3}) can now be established;
we fix $i$ and omit the superscript $(i)$ in the following discussion.
We consider the product of the following two matrices:
$$
\left(
\begin{array}{cc}
 \sigma I_{m} & -GBH^{\top} \\
- H B^{\top} G^{\top} & \sigma I_{n} \\
\end{array} \right) \,
\left(
\begin{array}{cc}
 \sigma I_m & 0 \\
 HB^{\top} G^{\top} & \sigma^{-1} I_n \\
\end{array}
\right),
$$
which is equal to
$$
\left(
\begin{array}{cc}
 \sigma^2 I_{m} -GBB^{\top} G^{\top}& -\sigma^{-1} GBH^{\top} \\
- \sigma H B^{\top} G^{\top} +\sigma HB^{\top} G^{\top} & I_{n} \\
\end{array} \right).
$$
Since the bottom-left block is $\mathbf{0}$, the determinant of this matrix
is $\det(\sigma^2 I_m - G B B^{\top} G)$.
Setting $C=BB^{\top}$ and $\tilde G = \left(g | G \right)$,
we have
$$
\sigma^2 I_m - GCG^{\top} =
\sigma^2 I_m - {\tilde G} 
\left(
\begin{array}{cc}
 0 & 0 \\ 
 0 & C \\
\end{array}
\right)
{\tilde G}^{\top}.
$$
Since ${\tilde G} {\tilde G}^{\top} = E$,
the determinant of the matrix above is equal to
$\sigma^2 \, \det(\sigma^2 I_{m-1} -C)$.
In summary, we  have obtained equalities of (\ref{eq:det2}) and (\ref{eq:det3}).

With regard to the expectation of the Euler characteristic number, exchanging the sum and the integral, we have
\[
{\rm E}[\chi(M_x)]
=
\sum_{i=1}^m \int dA p(A) \1(\sigma^{(i)} \geq x) \,\sgn\,
{\sigma^{(i)}}^{n-m} {\sigma^{(i)}}^2
\det
\left( 
{\sigma^{(i)}}^2 I_{m-1} - B^{(i)} {B^{(i)}}^{\top}
\right).
\]

To evaluate the expectation of the Euler characteristic number,
we require the Jacobian of 
(\ref{eq:ghB:coordinate}).
According to standard arguments in multivariate analysis (see, \eg, \cite[(3.19)]{takemura-kuriki-1999}),
we have
\[
  dA = \left| \det\left(\sigma^2 I_{m-1}-BB^{\top}\right)\right|\, 
      d\sigma G^{\top} dg H^{\top} dh dB.
\]
Subsequently, we have
\begin{equation*}
 {\rm E}[\chi(M_x)] \\
= \frac{1}{2} \sum_{i=1}^m \int_x^\infty \sigma^{n-m} d\sigma 
 \int_{B \in {\cal B}(i,\sigma^{(i)})} dB \int_{S^{m-1}} G^{\top} dg
 \int_{S^{n-1}} H^{\top} dh \ 
 \det\left(\sigma^2 I_{m-1}-B B^{\top}\right) p(A).
\end{equation*}
The factor $1/2$  is owing to the multiplicity of
$(g,h) \mapsto gh^{\top}$ being $2$.
Set ${\cal B}^{(i)} = {\cal B}(i,\sigma^{(i)})$.
For $i \neq j$, since ${\cal B}^{(i)} \cap {\cal B}^{(j)}$ and $\mathbb{R}^{(m-1)(n-1)} \setminus \sum_{i=1}^m {\cal B}^{(i)}$ are measure zero sets,
we may sum up integral domains for $B$ into one domain as
\[
 \sum_{i=1}^m \int_{B \in {\cal B}(i,\sigma^{(i)}) } \det(\sigma^2 I_{m-1}-B B^{\top}) p(A)
= \int_{B \in M(m-1,n-1) } \det(\sigma^2 I_{m-1}-B B^{\top}) p(A).
\]
Thus, we  have derived the conclusion.
\qedsymbol

The integral (\ref{EchiMx}) does not depend on  the choice of $G(g)$
nor $H(h)$. The reason is as follows.
The column vectors of the matrix $G=G(g)$  are of length 1 and 
are orthogonal to the vector $g$.
Let ${\tilde G}$ be a matrix with the same property.
In other words, we assume 
$(g,{\tilde G}) \in S\hspace*{-0.5mm}O(m)$.
Then there exists an $(m-1) \times (m-1)$ orthogonal matrix $P$
such that ${\tilde G}=G P$ and $|P|=1$ hold.
Taking the exterior product of elements of 
${\tilde G}^{\top} dg = P G^{\top} dg$, 
we have
$$ \wedge_{i=1}^m {\tilde g}_i^{\top} dg = |P| \, \wedge_{i=1}^m g_i^{\top} dg 
= \wedge_{i=1}^m g_i^{\top} dg.
$$
The case for $H$ can be shown analogously.

One of the most important examples is that $A$ has a Gaussian distribution $\mathcal{N}_{m\times n}(M,\Sigma\otimes I_n)$, 
 where $\otimes$ is the Kronecker product of matrices. In this case, we have
\begin{equation}
\label{EQ:gaussdis}
 p(A) dA = \frac{1}{(2\pi)^{mn/2}\det(\Sigma)^{n/2}}\exp\left\{-\frac{1}{2}\tr(A-M)^{\top}\Sigma^{-1}(A-M)\right\} dA.
\end{equation}
The largest singular value of $A$ is the square root of the largest eigenvalue of 
a non-central Wishart matrix $W_m(n,\Sigma,\Sigma^{-1}M M^{\top})$.
Substituting~\eqref{eq:ghB:coordinate} and~\eqref{EQ:gaussdis} into (\ref{EchiMx}), we have 
\begin{align}
 {\rm E}[\chi(M_x)]
 =& \frac{1}{2} \int_x^\infty \sigma^{n-m} d\sigma \int_{\mathbb{R}^{(m-1)\times(n-1)}} dB \int_{\S^{m-1}} G^{\top} dg \int_{\S^{n-1}} H^{\top} dh \det\left(\sigma^2 I_{m-1} - BB^{\top} \right) \nonumber \\
& \times\frac{1}{(2\pi)^{nm/2}\det(\Sigma)^{n/2}}\exp\left\{ -\frac{1}{2}\tr\left((\sigma h g^{\top} + H B^{\top} G^{\top} -M^{\top})\Sigma^{-1}(\sigma g h^{\top} + G B H^{\top} -M)\right) \right\}.
\label{eq:chi-for-normal}
\end{align}

In this expression, the number of parameters is $m(m+1)/2+mn$; therefore, it is over-parameterized.
Note that
\begin{equation}
\label{A}
 A = \Sigma^{1/2} V + M, \quad V=(v_{ij})_{m\times n}, \ \ v_{ij}\sim \mathcal{N}(0,1) \ \ {\rm i.i.d.}
\end{equation}
Let $\Sigma^{1/2}=P^{\top} D P$ be a spectral decomposition, where $D=\diag(d_i)$.
Then we have
$$ PA=DPV+PM. $$
Let $PM=NQ$ be a QR decomposition, where $N$ is $m\times n$ lower triangular matrix with nonnegative diagonal elements and $Q\in O(n)$.
Then $PAQ^{\top}=DV+N$.
Since the largest eigenvalues of $A$ and $PAQ^{\top}$ are the same, we can assume without loss of generality that $\Sigma$ is a diagonal matrix and $M$ is 
a lower triangle with nonnegative diagonal elements,~\ie,
\begin{equation}
 \Sigma^{-1} = \begin{pmatrix}
 s_1 &        & 0 \\
     & \ddots &   \\
 0   &        & s_m
 \end{pmatrix}, \ \ s_i>0, \quad
 M = \begin{pmatrix}
 m_{11} &        & 0      & 0      & \cdots & 0      \\
 \vdots & \ddots &        & \vdots &        & \vdots \\
 m_{mn} & \cdots & m_{mm} & 0      & \cdots & 0
 \end{pmatrix}, \ \ m_{ii}\ge 0.
\label{parameters}
\end{equation}
When $\Sigma$ has multiple roots, \ie,
\begin{equation}
\label{parameters1}
 \Sigma^{-1} = \begin{pmatrix}
 s_1 I_{n_1} &        & 0 \\
             & \ddots &   \\
 0           &        & s_r I_{n_r}
\end{pmatrix}, \quad \sum_{i=1}^r n_i=m,
\end{equation}
by multiplying $\diag(P_1,\ldots,P_r)\in O(n_1)\times\cdots\times O(n_r)$ and its transpose from the left and right, we can assume
\begin{equation}
\label{parameters2}
 M = \begin{pmatrix}
 m_1 I_{n_1} &             &        &                     &             & 0      & \cdots & 0 \\
 M_{21}      & m_2 I_{n_2} &        &                     &             & 0      & \cdots & 0 \\
 \vdots      & \vdots      & \ddots &                     &             & \vdots &        & \vdots \\
 M_{r-1,1}   & M_{r-1,2}   &        & m_{r-1} I_{n_{r-1}} &             & 0      & \cdots & 0 \\
 M_{r1}      & M_{r2}      & \cdots & M_{r,r-1}           & m_r I_{n_r} & 0      & \cdots & 0
\end{pmatrix}, \ \ m_i\ge 0, \ \ M_{ij}\in\mathbb{R}^{n_j\times n_i}.
\end{equation}
Therefore, our problem can be formalized as follows: Evaluate (\ref{eq:chi-for-normal}) with 
parameters (\ref{parameters}) (or (\ref{parameters1}) and (\ref{parameters2})).

In the following sections, we will evaluate the integral representation
of the expectation of the Euler characteristic
given in Theorem \ref{th:expectation_of_euler_characteristic}
for some interesting special cases. 
We can obtain approximate values of the probability
of the largest eigenvalue of random matrices by virtue of them.
The Euler characteristic heuristic is
\[
 \Pr\left(\max_{g\in S^{m-1}, h\in S^{n-1}} g^{\top} A h \ge x\right)
 = \Pr\left(\max_{U\in M}f(U)\ge x\right)
 \approx {\rm E}\left[\chi(M_x)\right].
\]
The condition that $f(U)$ is a Morse function with probability one holds if $A$ has  $m$ distinct and non-zero singular values with probability one.

\section{The case of $m=n=2$}
\label{sec:case}

We derive Theorem \ref{th:expectation_of_euler_characteristic}
in the special case of $m=n=2$ by taking explicit coordinates.
This derivation motivates the proof for the general case
discussed in the previous section. The case $m=n=2$ is studied numerically in the last section with 
the holonomic gradient method (HGM).

Fix two unit vectors 
\[
g=\left(\cos\theta,\sin\theta\right)^{\top},h=\left(\cos\phi,\sin\phi\right)^{\top}\in S^{1}, \ \ \text{ for } \ \ 0\leq \theta, \phi<2\pi. 
\]
Define 
\[
G=\left(\cos\left(\theta+\frac{\pi}{2}\right),\sin\left(\theta+\frac{\pi}{2}\right)\right)^{\top}=\left(-\sin\theta,\cos\theta\right)^{\top},
\]
which satisfies 
\[
\left(g,G\right)=\left(\begin{matrix}\cos\theta & -\sin\theta\\
\sin\theta & \cos\theta
\end{matrix}\right)\in \mathit{SO}(2).
\]
Similarly, we define $H=\left(\cos\left(\phi+\frac{\pi}{2}\right),\sin\left(\phi+\frac{\pi}{2}\right)\right)^{\top}=\left(-\sin\phi,\cos\phi\right)^{\top}$.
Here,  in case the sum is greater than $2\pi$, both $\theta+\pi/2$ and $\phi+\pi/2$ should
be treated as mod $2\pi$.
Now, any $2\times2$ matrix, say $A$, can be recovered by 
\[
A=\sigma gh^{\top}+bGH^{\top}
\]
with $4$ variables $\left(\sigma,\theta,\phi,b\right)$. 
We may further assume that $\sigma \in \mathbb{R}_{\geq0}$, $b\in\mathbb{R}$, and $\phi,\theta\in[0,2\pi)$.

Fix $\sigma_{0},b_{0},\theta_{0}, \phi_{0}$ and let
\[
A_{0} =\sigma_{0}g\left(\theta_{0}\right)h\left(\phi_{0}\right)^{\top}+b_{0}G\left(\theta_{0}\right)H\left(\phi_{0}\right)^{\top}.
\]
By  letting $\sigma,b$ vary in $\mathbb{R}$ and $\phi,\theta$ vary
in $[0,2\pi)$, we recover $A_{0}$ four times:
\[
\left\{
\begin{aligned}
A_{0} & =\sigma_{0}g\left(\theta_{0}\right)h\left(\phi_{0}\right)^{\top}+b_{0}G\left(\theta_{0}\right)H\left(\phi_{0}\right)^{\top};\\
A_{0} & =\sigma_{0}g\left(-\theta_{0}\right)h\left(-\phi_{0}\right)^{\top}+b_{0}G\left(-\theta_{0}\right)H\left(-\phi_{0}\right)^{\top};\\
A_{0} & =b_{0}g\left(\theta_{0}+\frac{\pi}{2}\right)h\left(\phi_{0}+\frac{\pi}{2}\right)^{\top}+\sigma_{0}G\left(\theta_{0}+\frac{\pi}{2}
\right)H\left(\phi_{0}+\frac{\pi}{2}\right)^{\top};\\
A_{0} & =b_{0}g\left(-\theta_{0}+\frac{\pi}{2}\right)h\left(-\phi_{0}+\frac{\pi}{2}\right)^{\top}+\sigma_{0}G\left(-\theta_{0}
+\frac{\pi}{2}\right)H\left(-\phi_{0}+\frac{\pi}{2}\right)^{\top}.
\end{aligned}
\right.
\]

\begin{itemize}
\item Here, for the first two cases, it is easily seen from the symmetry of the manifold
$M$ (shown below) that $\left(h,g\right)$ is equivalent to $\left(-h,-g\right)$.
\item The second symmetry is given by $\left(\sigma',b'\right)=\left(b_{0},\sigma_{0}\right)$,
\ie, interchanging $\sigma$ and $b$. Note that $G\left(\theta\right)=g\left(\theta+\frac{\pi}{2}\right)$
and $H\left(\phi\right)=h\left(\phi+\pi/2\right)$. Thus, there also exists 
$$ \left(\theta',\phi'\right)=\left(\theta_{0}+\frac{\pi}{2},\phi_{0}+\frac{\pi}{2}\right) $$
recovering $A_{0}$. 
\end{itemize}
Therefore, to recover $A$, we can always assume
that $\sigma\geq b$ and let $\theta,\phi\in[0,2\pi)$.
See Lemma \ref{lemma:one-to-one-2017-12-08} for a general claim.

Next, we consider the manifold 
\[
M =\left\{ ts^{\top} \mid s=\left(\cos\alpha,\sin\alpha\right),t=\left(\cos\beta,\sin\beta\right)\in S^{1},0\leq\alpha,\beta<2\pi\right\} 
\]
and the function $f$ on $M$ such that 
\[
f\left(ts^{\top}\right) =s^{\top}At=s^{\top}\left(\sigma gh^{\top}+bGH^{\top}\right)t.
\]
Apparently, $A$ only has two pairs of eigenvectors, which can be verified
by the following computations:
\[
\left\{
\begin{aligned}
 Ah=& \sigma gh^{\top}h+bGH^{\top}h =\sigma g;\\
 g^{\top}A=& \sigma g^{\top}gh^{\top}+bg^{\top}GH^{\top} =\sigma h^{\top};\\
 AH=& \sigma gh^{\top}H+bGH^{\top}H =bG;\\
 G^{\top}A=& \sigma G^{\top}gh^{\top}+bG^{\top}GH^{\top} =bH^{\top}.
\end{aligned}
\right.
\]
The function $f$ has two critical points on $M$, which are at 
\begin{itemize}
\item the point $P=hg^{\top}\in M\Leftrightarrow\left(\alpha,\beta\right)=\left(\theta,\phi\right)$ 
\item or $Q=HG^{\top}\in M\Leftrightarrow\left(\alpha,\beta\right)=\left(\theta+\pi/2,\phi+\pi/2\right)$. 
\end{itemize}
Further computation indicates the following four facts:
\begin{enumerate}
\item[(i)] $f\left(P\right)=g^{\top}Ah=\sigma$ and $f\left(Q\right)=G^{\top}AH=b$;
\item[(ii)] From 
\begin{align*}
 \mathrm{Hess}f & = \left(\begin{matrix}\frac{\partial^{2}}{\partial\alpha^{2}}f & \frac{\partial^{2}}{\partial\alpha\partial\beta}f\\
\frac{\partial^{2}}{\partial\beta\partial\alpha}f & \frac{\partial^{2}}{\partial\beta^{2}}f
\end{matrix}\right) 
 = \frac{1}{2}
 \left(\begin{matrix}
   b - \sigma & b + \sigma \\
   b - \sigma & - b - \sigma
 \end{matrix}\right) 
 \left(\begin{matrix}
   \cos(\alpha + \beta - \theta - \phi) & \cos(\alpha + \beta - \theta - \phi) \\
   -\cos(\alpha - \beta - \theta + \phi) & \cos(\alpha - \beta - \theta + \phi)
 \end{matrix}\right),
\end{align*}
it follows that $\det\left(\mathrm{Hess}_{P}f\right)=\sigma^{2}-b^{2}$
and $\det\left(\mathrm{Hess}_{Q}f\right)=b^{2}-\sigma^{2}$. Therefore,
we see
\begin{enumerate}
\item if $x>\sigma\geq b$, then $M_{x}$ does not contain any critical points,
so $\chi\left(M_{x}\right)=0$;
\item if $x<b\leq\sigma$, then $M_{x}$ contains both critical points, and thus
\[
\chi\left(M_{x}\right)=\sgn\left(\sigma^{2}-b^{2}\right)+\sgn\left(b^{2}-\sigma^{2}\right)=0;
\]
\item the only nontrivial case is $\sigma\geq x\geq b$, then 
\[
\chi\left(M_{x}\right) = \1\left(\sigma\geq x\geq b\right)\sgn\left(\sigma^{2}-b^{2}\right).
\]
\end{enumerate}
\item[(iii)] Since
\[
 A=\sigma gh^{\top}+bGH^{\top}
=\left(\begin{array}{cc}
 b\sin\theta\sin\phi+\sigma\cos\theta\cos\phi & \sigma\cos\theta\sin\phi-b\sin\theta\cos\phi\\
 \sigma\sin\theta\cos\phi-b\cos\theta\sin\phi & b\cos\theta\cos\phi+\sigma\sin\theta\sin\phi
\end{array}\right),
\]
we have 
\begin{align*}
\left(dA\right)
&= db\sin\theta\sin\phi+\sigma\cos\theta\cos\phi\wedge d\left(\sigma\cos\theta\sin\phi-b
\sin\theta\cos\phi\right) \\ 
& \ \ \ \wedge d\left(\sigma\sin\theta\cos\phi-b\cos\theta\sin\phi\right)\wedge d\left(b\cos\theta\cos\phi
+\sigma\sin\theta\sin\phi\right)\\
&= \left(b^{2}-\sigma^{2}\right) d\sigma db d\theta d\phi,
\end{align*}
where $\wedge$ is the exterior product for vectors.
\item[(iv)] Let $M=\left(\begin{matrix}m_{11} & 0\\
m_{21} & m_{22}
\end{matrix}\right)$ and $\Sigma=\left(\begin{matrix}1/s_{1} & 0\\
0 & 1/s_{2}
\end{matrix}\right)$ such that 
\[
 A=\sqrt{\Sigma}V+M,\text{ where }V=\left(v_{ij}\right),\ v_{ij}\sim\mathcal{N}\left(0,1\right)\text{ i.\,i.\,d.}
\]
Then
\[
 p\left(A\right)=\frac{s_{1}s_{2}}{\left(2\pi\right)^{2}}e^{-\frac{R}{2}},
\]
where 
\begin{align*}
 R=
 & s_{1}\left(b\sin\theta\sin\phi+\sigma\cos\theta\cos\phi-m_{11}\right)^{2}+s_{2}\left(\sigma\sin\theta\cos\phi-b\cos\theta\sin\phi-m_{21}
\right)^{2}\\ 
 & +s_{1}\left(\sigma\cos\theta\sin\phi-b\sin\theta\cos\phi\right)^{2}
 +s_{2}\left(b\cos\theta\cos\phi+\sigma\sin\theta\sin\phi-m_{22}\right)^{2}.
\end{align*}
\end{enumerate}
Hence, we have 
\begin{align*}
{\rm E}\left[\chi\left(M_{x}\right)\right]
=& \frac{1}{2}\int_{-\infty}^{\infty} d\sigma\int_{-\infty}^{\infty}
 db\int_{0}^{2\pi} d\theta\int_{0}^{2\pi} d\phi\left(\1\left(\sigma\geq x\geq b\right)
\sgn\left(\sigma^{2}-b^{2}\right)\right) \left|\left(b^{2}
 -\sigma^{2}\right)\right|
 \frac{s_{1}s_{2}}{\left(2\pi\right)^{2}}e^{-\frac{R}{2}} \\
=& \frac{1}{2}\int_{x}^{\infty} d\sigma\int_{-\infty}^{x} db\int_{0}^{2\pi} d\theta\int_{0}^{2\pi} d\phi
 \left(\sigma^{2}-b^{2}\right)\frac{s_{1}s_{2}}{\left(2\pi\right)^{2}}e^{-\frac{R}{2}}.
\end{align*}
Note that we have
$\int_{-\infty}^\infty db \ldots = \int_{-\infty}^x db \ldots$
by the anti-symmetry of $\sigma$ and $b$ in this case.
In other words, integrals over 
$\sigma > x >0, b>x, \sigma > b$
and 
$\sigma > x >0, b>x, \sigma < b$
are canceled. Thus, we have
\begin{equation}
\label{EQ:mn2}
 {\rm E}[\chi(M_x)] = F(s_1,s_2,m_{11},m_{21},m_{22};x) 
= \frac{1}{2}\int_x^\infty d\sigma \int_{-\infty}^\infty db \int_0^{2\pi} d\theta \int_0^{2\pi} d\phi (\sigma^2-b^2) 
\frac{s_1 s_2}{(2\pi)^2} \exp\left(-\frac{1}{2}R \right). 
\end{equation}
In summary, we have obtained Theorem \ref{th:expectation_of_euler_characteristic}
in the case that $A$ has a Gaussian distribution.

A numerical example is given below.
\begin{example}
\label{EXM1}\rm
We evaluate~\eqref{EQ:mn2} with parameters 
$$ s_1=2,\ \ s_2= m_{11}=1,\ \ m_{21}=-1,\ \ m_{22}=1 $$
and derive Table~\ref{TABLE:1}.
\begin{table}[h]
\caption{Euler characteristic versus Monte Carlo simulation for the evaluation of~\eqref{EQ:mn2}}
\begin{center}
\begin{tabular}{ccccccccc}
\hline
$x$ & 0 & 1 & 2 & 3 & 4 & 5 \\
\hline
${\rm E}[\chi(M_x)]$ & $-5.9\times 10^{-8}$ & 0.74 & 0.56 & 0.14 & 0.014 & 0.00058 \\
%
$\Pr(\sigma>x)$ & 1. & 0.95 & 0.57 & 0.14 & 0.014 & 0.00058 \\
\hline
\end{tabular}
\end{center}
\label{TABLE:1}
\end{table}
Here, the probability $\Pr(\sigma>x)$ is estimated by a Monte Carlo simulation with 10,000,000 iterations, and the expectation of the Euler characteristic is evaluated by a numerical
integration function {\tt NIntegrate} on Mathematica.
As expected, ${\rm E}[\chi(M_x)]\approx \Pr(\sigma>x)$ when $x$ is large.
\end{example}

\section{Computer algebra and the expectation for small $m$ and $n$}
\label{sec:computer}

In this section, we study the non-central case
$M \neq 0$ with the help of computer algebra.
When $m = n = 2$, we can perform the holonomic gradient method (HGM)~\cite{Hashiguchi2013}
to evaluate the integral~\eqref{EchiMx}.


In Section~\ref{sec:case}, we derive an integral formula~\eqref{EQ:mn2}  for the case $m = n = 2$. For~\eqref{EQ:mn2}, 
we set 
$$
 \sin \theta = \frac{2 s}{1 + s^2},\ \ \cos \theta = \frac{1 - s^2}{1 + s^2}, 
 \ \ \sin \phi = \frac{2 t}{1 + t^2},\ \ \cos \phi = \frac{1 - t^2}{1 + t^2}.
$$ 
Then we have 
\begin{equation} 
 {\rm E}[\chi(M_x)] = F(s_1,s_2,m_{11},m_{21},m_{22};x) 
= \frac{1}{2 \pi^2}\int_x^\infty d\sigma \int_{-\infty}^\infty db 
\int_{-\infty}^{\infty} ds \int_{-\infty}^{\infty} dt \frac{s_1 s_2 (\sigma^2-b^2)}{(1 + s^2) (1 + t^2)} 
\exp\left(-\frac{1}{2} \tilde{R} \right),
\label{EQ:secondeuler}
\end{equation} 
where $\tilde{R}$ is a rational function in $\sigma, b, s, t$. 
Since the integrand is a holonomic function in $\sigma, b, s, t$, 
we can apply the creative telescoping method~\cite{Zeilberger1991} to derive holonomic systems for the integrals. 
This is straightforward for the inner single integral of ${\rm E}[\chi(M_x)]$ by 
the classic methods~\cite{Christoph2009} 
(such as Zeilberger's algorithm, Takayama's algorithm and Chyzak's algorithm). 
Below is an example.

\begin{example}
\label{EX:firstintegral} \rm
Consider the inner single integral of~\eqref{EQ:secondeuler}:
\[
 f_1(\sigma, b, s)
 = \int_{-\infty}^{\infty} \frac{s_1 s_2 (\sigma^2-b^2)}{(1 + s^2) (1 + t^2)} 
 \exp\left(-\frac{1}{2} \tilde{R} \right) dt
 = \int_{-\infty}^{\infty} f_0 \cdot dt,
\]
where $\tilde{R}$ is a rational function in $\sigma, b, s, t$. 
 Since $f_0$ is a holonomic function, we can compute a holonomic system satisfied by $f_0$
using the Mathematica package \texttt{HolonomicFunctions}~\cite{Christoph2010}.
Using the holonomic system satisfied by $f_0$ and Chyzak's algorithm~\cite{Chyzak2000}, 
we can then derive a holonomic system of $f_1$, which is of holonomic rank $2$. 
The detailed calculations can be found in the supplementary material~\cite{ElectronicYZ}.
%
%
%
%
%
\end{example}
In the aforementioned example, we use Chyzak's algorithm to derive a holonomic system of the inner single integral of~${\rm E}[\chi(M_x)]$. 
This can be done within 5 seconds on a Linux computer with 15.10 GB RAM. 
However, experiments show that it is not efficient enough to derive a holonomic system for the inner double 
integral in the same way within reasonable computational time because of the complexity of this algorithm. 
To speed up the computation, 
we intend is to utilize the Stafford theorem~\cite{Hillebrand2001, Anton2004} empirically. 
Let us first recall the theorem. 
Assume that $\mathbb{K}$ is a field of characteristic $0$ and~$n$ is a positive integer. 
Let $R_n = \mathbb{K}(x_1, \ldots, x_n)[\partial_1, \ldots, \partial_n]$ 
and $D_n = \mathbb{K}[x_1, \ldots, x_n][\partial_1, \ldots, \partial_n]$ be 
the ring of differential operators with rational coefficients and the Weyl algebra in $n$ variables, respectively. 

\begin{theorem}
\label{THM:Stafford} 
Every left ideal in $R_n$ or $D_n$ can be generated by two elements. 
\end{theorem}

Assume that $I$ is a left ideal in $R_n$ or $D_n$. 
We observe from experiments that for any two random operators $a, b \in I$, 
it is of high probability that $I = \langle a, b \rangle$. 
This suggests the following heuristic method for computing a holonomic system 
for the inner double integral of ${\rm E}[\chi(M_x)]$. 
As a matter of notation, we set 
$$
 T_{n - 1} = \{ \partial_1^{i_1} \partial_2^{i_2} \cdots \partial_{n - 1}^{i_{n -1}} 
\mid (i_1, \ldots, qi_{n - 1}) \in \mathbb{N}^{n - 1} \}.
$$

Recall that a D-finite system~\cite{Yi2017} in $R_n$ is 
a finite set of generators of a zero-dimensional ideal in $R_n$. 
The relation between D-finite systems and holonomic systems is illustrated in~\cite[Section 6.9]{Hibi2013}. 
For the application of the HGM, 
D-finite systems are alternative to holonomic systems. 
Here, we use D-finite systems because they are more efficient for computation. 


\begin{heuristic}
\label{HEU:stafford} \rm
Given a D-finite system $G$ in $R_n$, 
compute another D-finite system $G_1$ in~$R_{n - 1}$ 
such that
$$ G_1 \subset \left( R_n \cdot G + \partial_n R_n \right) \cap R_{n - 1}. $$
\begin{enumerate}
 \item[(i)] Choose two finite support set $S_1, S_2 \in T_{n - 1}$. 
 \item[(ii)] Using the polynomial ansatz method~\cite[Section 3.4]{Christoph2009}, 
 check whether there exist telescopers $P_1, P_2 \in R_{n - 1}$ of $G$ with support sets $S_1, S_2$ or not. 
 If $P_1$ and $P_2$ exist, then go to the next step. Otherwise, go to step 1. 
 \item[(iii)] Compute the Gr\"{o}bner basis $G_1$ of $\{P_1, P_2 \}$ with respect to a term order~\cite{Cox2015} in $T_{n - 1}$. 
 If $G_1$ is D-finite, then output $G_1$. Otherwise, go to step 1. 
\end{enumerate}
\end{heuristic}

In the aforementioned heuristic method, we need to find two finite support set $S_1, S_2 \in T_{n - 1}$ 
through trial and error so that the computation terminates and finishes in reasonable time.
Next, we demonstrate its application to derive a D-finite system for the inner double 
integral of ${\rm E}[\chi(M_x)]$.

\begin{example}
\label{EX:stafford heuristic} \rm
Consider the inner double integral of~\eqref{EQ:secondeuler}:
\begin{equation}
\label{EQ:secondintegral}
 f_2(\sigma, b) = \int_{-\infty}^{\infty} f_1(\sigma, b, s) ds
\end{equation}
where $f_1(\sigma, b, s)$ is defined in Example~\ref{EX:firstintegral}. 

Let $G$ be a D-finite system of $f_1$, which is derived from Example~\ref{EX:firstintegral}. 
Using $G$ and the polynomial ansatz method, we find two non-zero annihilators $P_1$ and $P_2$ for $f_2$ with support sets $S_1$ and $S_2$, 
respectively, where 
\[
S_1 = \{1, \partial_b, \partial_{\sigma}, \partial_{b}^2, \partial_b \partial_{\sigma}, \partial_{\sigma}^2, \partial_{\sigma}^3\}, \quad \quad S_2 = S_1 \cup \{ \partial_b^2 \partial_{\sigma}, \partial_b \partial_{\sigma}^2, \partial_b^3\}.
\]
Then we compute the Gr\"{o}bner basis $G_1$ of $\{P_1, P_2 \}$ in $\mathbb{Q}(b, \sigma)[\partial_b, \partial_{\sigma}]$ with respect to 
a total degree lexicographic order. We find that $G_1$ is a D-finite system of holonomic rank $6$. 
The details of the calculation can be found in~\cite{ElectronicYZ}.
%
%
%
%
%
%
%
%
%
%
%
%
%
%
\end{example}

In the aforementioned example, we specify the parameters in the integrand as those in Example~\ref{EXM1}. 
Using Heuristic~\ref{HEU:stafford}, we can further compute a holonomic system for 
the inner double integral of~${\rm E}[\chi(M_x)]$ without specifying those parameters (pars). 
This is significantly more efficient than Chyzak's algorithm. 
Table~\ref{Table:2} compares Chyzak's algorithm (chyzak) and Heuristic~\ref{HEU:stafford} (heuristic) 
in terms of computation time (s).

\begin{table}[h]
\caption{Chyzak's algorithm versus Heuristic~\ref{HEU:stafford} in deriving holonomic systems of~\eqref{EQ:secondintegral}}
\begin{center}
\begin{tabular}{ c c c c c c c }
\hline
\# pars  & $0$     & $1$                   & $2$    & $3$    & $4$          & $5$ \\
\hline
Chyzak   & $976$   & $9.8323\times 10^{4}$ & -      & -      & -            & - \\
Heuristic& $43.49$ & $394.4$               & $8527$ & $4.3957\times 10^{5}$ & - & $1.5519\times 10^{6}$\\
\hline
\end{tabular}
\end{center}
\label{Table:2}
\end{table}

Next, we use Heuristic~\ref{HEU:stafford} to derive a D-finite system of the inner triple integral of ${\rm E}[\chi(M_x)]$ 
and then numerically solve the corresponding ordinary differential equation. 
Finally, we use numerical integration to evaluate ${\rm E}[\chi(M_x)]$.

\begin{example}
\label{EX:hgm} \rm
Consider
\begin{equation}
\label{EQ:eulerexpectation}
 {\rm E}[\chi(M_x)] = \frac{1}{2 \pi^2}\int_x^\infty d\sigma \int_{-\infty}^\infty db f_2(\sigma, b),
\end{equation}
where $f_2(\sigma, b)$ is specified in~\eqref{EQ:secondintegral} with parameters 
\[
 s_1=2, s_2= m_{11}=1, m_{21}=-1, m_{22}=1.
\] 

By Example~\ref{EX:stafford heuristic}, we have derived a D-finite system for~$f_2$. 
Using Heuristic~\ref{HEU:stafford}, 
we derive a D-finite system for the inner first integral $f_3$ of~\eqref{EQ:eulerexpectation} 
of the following form:
\[
 P = c_{10} \cdot \partial_{\sigma}^{10} + c_{9} \cdot \partial_{\sigma}^9 + \cdots + c_{0},
\]
where $c_i \in \mathbb{Q}[\sigma], i \in \{ 0, \ldots, 10\}$.
%
%
%
%
%
Now, we first numerically solve the ordinary differential equation $P(f_3) = 0$ to evaluate $f_3$, 
and then evaluate ${\rm E}[\chi(M_x)]$ by using numerical integration. 
Table~\ref{TABLE:3} are the corresponding numerical results,
\begin{table}[h]
\caption{Holonomic gradient method (HGM) versus Monte Carlo simulation in evaluating ${\rm E}[\chi(M_x)]$}
\begin{center}\small
\begin{tabular}{cccccccc}
\hline
$x$ & $1$        & $2$        & $3$        & $4$         & $5$           & $6$ \\
\hline
HGM & $0.745835$ & $0.567729$ & $0.144879$ & $0.0146728$ & $0.000582526$ & $8.79942\times 10^{-6}$ \\
mc  & $0.745802$ & $0.567623$ & $0.144986$ & $0.0146901$ & $0.0005933$    & $9.6\times 10^{-6}$ \\
\hline
\end{tabular}
\end{center}
\label{TABLE:3}
\end{table}
 where mc represents the Monte Carlo simulation of ${\rm E}[\chi(M_x)]$ by the following formula with 10,000,000 iterations:
\[
 {\rm E}[\chi(M_x)] \approx \frac{\sum_{i = 1}^{m} \chi(M_{x, i})}{m},
\]
with 
$$ \chi(M_{x, i}) = \1(\sigma_i \ge x) (\sigma_i^2 - b_i^2) + \1(b_i \ge x) (b_i^2 - \sigma_i^2), $$ 
where $\sigma_i$ and $b_i$ are singular values of $M_{x, i}$, $i \in \{ 1, \ldots, m\}$.

As expected, the results of the HGM are approximate to those of the Monte Carlo simulation. 
The detailed computation can be found in~\cite{ElectronicYZ}.
\end{example}

The evaluations of ${\rm E}[\chi(M_x)]$ in the above example are also approximate of those given in Example~\ref{EXM1}. 
The source codes for this section and a demo notebook are freely available as part of the supplementary
electronic material~\cite{ElectronicYZ}.




\begin{example} \rm
We consider the evaluation of~\eqref{EQ:secondeuler} with parameters
$$ m_{11}=1,\ \ m_{21}=2,\ \ m_{22}=3,\ \ s_1=10^3,\ \ s_2=10^2. $$
It is difficult to evaluate~\eqref{EQ:secondeuler} 
for the relatively large parameters $s_i$ by numerical integration (even with the Monte Carlo integration).
Thus, we take a different approach. 
Using Heuristic~\ref{HEU:stafford}, we can compute a linear ordinary differential equation (ODE) for~\eqref{EQ:secondeuler} of rank $11$ with respect to the independent variable $x$.   
Then we construct series solutions for this differential
equation and use them to extrapolate results by simulations.

Although this extrapolation method is well known,
we explain it in a subtle form with application in our evaluation problem.
Consider an ODE with coefficients in $\mathbb{Q}(x)$ of rank~$r$.
Let $c \in \mathbb{Q}$ be a point in the $x$-space and 
we take $r$ increasing numbers $y_j \in \mathbb{Q}$, where $j \in \{0, \ldots, r-1\}$.
We construct a series solution $f_i(x)$
as a series in $x-(c+y_i)$.
We may further assume that $c+y_i$ is not a singular point of the ODE for each $i$.
The initial value vector may be taken suitably so that the series
is determined uniquely over $\mathbb{Q}$.

We assume that the vector $(f_i(x))$ converges in a segment $I$ containing all $c+y_i$'s
and that it is a basis of the solution space.
Once we construct such a basis of series solutions,
we can construct the solution $f(x)$ that takes values
$b_j$ at $x=p_j \in \mathbb{Q}\cap I$, $j \in \{0, \ldots, r-1\}$.
To be specific, we set 
$$ 
 f(x) = \sum_{i=0}^{r-1} t_i f_i(x)
$$
with unknown coefficients $t_i$'s. 
Then we have
$$ 
 f(p_j) = \sum_{i=0}^{r-1} t_i f_i(p_j), \quad j \in \{0, \ldots, r-1\}.
$$
The unknown coefficients $t_i$'s can be determined by solving the system of linear equations
\begin{equation}
\label{eq:lineq20190124}
 b_j = \sum_{i=0}^{r-1} t_i f_i(p_j)
\end{equation}
We call $f$ the extrapolation function by series solutions of ODE.
We call $b_j$ the reference value of $f$ at reference point $p_j$.

Let us now come back to our example. 
The linear ODE for~\eqref{EQ:secondeuler} has rank $r = 11$.
We set $c=370/100-1/100$ and 
the $y_j$'s as $[0, 1/100, \ldots, 10/100]$.
Then we have 
$$ c+y_0=3.69,\ \ c+y_1=3.70,\ \ \ldots,\ \ c+y_{10}=3.79. $$
We construct an approximate series solution $f_i(x)$
by taking $20000$ terms with rational arithmetic.

We set the reference points $p_j = 38/100+j/1000$,
$p_0=3.8, \ldots, p_{10}=3.81$
and construct a matrix related to~\eqref{eq:lineq20190124}.
Numbers in the matrix are translated to approximate rational numbers
to avoid the instability problem of solving linear equations 
(\ref{eq:lineq20190124})
with floating point numbers.

We assume that the expectation of the Euler characteristic of $M_x$
is almost equal to the probability $\Pr(\ell_1>x)$ that the first eigenvalue
is larger than $x$.
In fact, we have the Euler expectation ${\rm E}[\chi(M_x)] = \Pr(\ell_1 > x) - \Pr(\ell_2 >x)$
in this case, where $\ell_i$ is the $i$-th eigenvalue.
We have $\Pr(\ell_2>3.8)=0$ by the Monte-Carlo simulation with $1,000,000$ tries.
Then we may suppose that 
reference values $f(p_j)$ are estimated by Monte-Carlo simulation 
for $\Pr(\ell_1>x)$.
We construct a solution $f(x)$ with these reference values.
Evaluation of $f(x)$ is done with big floats.

\begin{table}[h]
\begin{center}
\caption{Numerical evaluation by extrapolation series versus Monte Carlo simulation for ${\rm E}[\chi(M_x)]$}
\label{table:values2}
\begin{tabular}{lll}
\hline
 $x$    &  $f(x)$  & simulation \\ \hline 
 3.8133 & 0.051146 & 0.051176 \\  
 3.8166 & 0.047517 & 0.047695 \\  
 3.82   & 0.044120 & 0.044515 \\   
 \hline
\end{tabular}
\end{center}
\end{table}
Table~\ref{table:values2} represents the values of the extrapolation
function $f(x)$ 
obtained by the above method with the big floats of 380 digits
and that by simulation with $1,000,000$ samples.
One simulation takes approximately $573$s by using the \textsf{R} package~\texttt{ mnormt} on a machine with Intel Xeon CPU(2.70GHz) and 256G memory.

The solid line in Fig.~\ref{fig:far} is obtained by
this extrapolation function.
The line goes to a big value at $x=3.866$ because this $x$ is out of 
the domain of convergence of this approximate series.
Dots are values obtained by simulation and those on the thick solid line 
are values used as reference values to obtain the extrapolation function.

\begin{figure}[hbt!]
\begin{center}
\includegraphics[width=10cm]{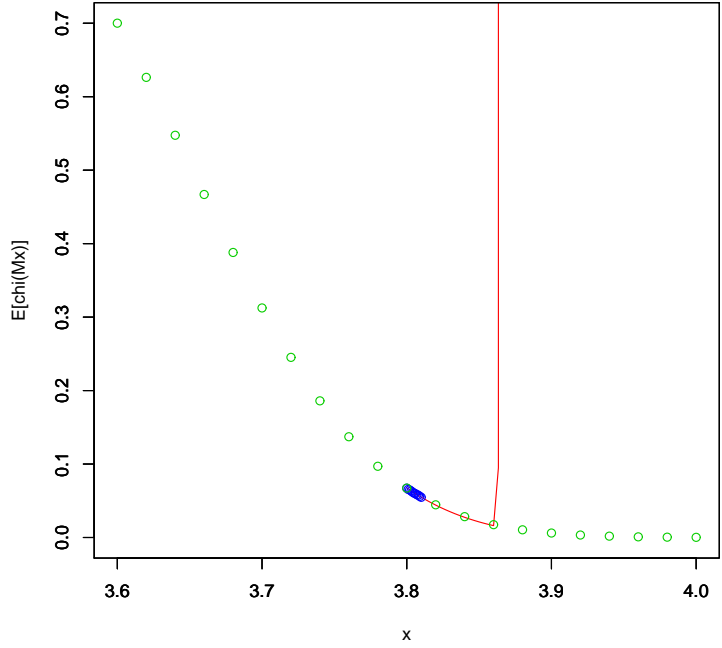}
\end{center}
\caption{Extrapolation function with 20000 terms.
Solid line is the extrapolation function,
which diverges when $x>3.8633$.
Dots are values from Monte Carlo simulation.
}
\label{fig:far}
\end{figure}

We obtain the series $f_i$ with $20,000$ terms in $5661$s by using Risa/Asir on a machine with Intel Xeon CPU(2.70GHz) and 256G memory.
The time to evaluate the extrapolation function at $61$ points is $14.03$s.
On the other hand, if we want to obtain simulation values at 61 points, we need about
$573 \times 61 = 34953s$.
Thus, our extrapolation method is of advantage in evaluating 
the function ${\rm E}[\chi(M_x)]$ for many $x$.
\end{example}

\appendix

\section{Proof of Lemma~\ref{lem:m=2}} \label{SEC:lemma2}

Recall that we are dealing with the Wishart matrix $W=A A^{\top}$ with $A$ given in (\ref{A}).
For an $m\times 1$ unit vector $p$,
$p^{\top} W p$ is distributed as
\[
 p^{\top}\Sigma p \cdot \chi^2\left(n;\frac{p^{\top} M p}{p^{\top}\Sigma p}\right), \qquad M=\Sigma\Omega,
\]
where
$c\cdot\chi^2(n;\delta^2)$ represents the distribution of $c$ times a non-central chi-square random variable with $n$ degrees of freedom and non-central parameter $\delta^2$.

We consider the case $m=2$.
From the characterization of the largest and smallest eigenvalues, we have
$\lambda_1(A)\ge p_1^{\top} A p_1$ and $\lambda_2(A)\le p_2^{\top} A p_2$, where $p_1$ and $p_2$ are arbitrary $2\times 1$  unit vectors.

(i) Suppose that $\Sigma$ has two distinct eigenvalues.
Set $p_1$ and $p_2$ to be two eigenvectors of $\Sigma$ corresponding to the eigenvalues $\lambda_1(\Sigma)$ and $\lambda_2(\Sigma)$.
Then
\begin{equation}
\label{ratio}
 \frac{\Pr(\lambda_2(W)\ge x)}{\Pr(\lambda_1(W)\ge x)}
 \le \frac{\Pr\bigl(\lambda_2(\Sigma) \chi^2(n;\delta_2^2)\ge x\bigr)}{\Pr\bigl(\lambda_1(\Sigma) \chi^2(n;\delta_1^2)\ge x\bigr)},
\end{equation}
where $\delta^2_i=p_i^{\top} M p_i/\lambda_i(\Sigma)$, $i=1,2$.
Note that $\delta_i$ can be zero. 

The tail behaviors of the central and non-central chi-square distributions were investigated by Beran \cite{beran-1975}.
From (2.9) and (3.3) of \cite{beran-1975}, combined with the asymptotics for the modified Bessel function of the first kind $I_\nu(x) \sim e^x/\sqrt{x}$ as $x\to\infty$, we have
\[
 \Pr\Bigl(\chi^2(n;b^2)\ge x\Bigr) \asymp
 \begin{cases}
  x^{(n-2)/2} e^{-x/2} & (b=0), \\
  x^{(n-3)/4} e^{-x/2+b \sqrt{x}} & (b>0),
 \end{cases}
\]
as $x\to\infty$.
In either case whether $\delta_i^2$ is zero or not, the right-hand side of (\ref{ratio}) goes to zero as $x$ goes to infinity.

(ii) Suppose that $\Sigma=\sigma I_2$ and $M=\Sigma\Omega$ has two distinct eigenvalues.
Set $p_1$ and $p_2$ to be two eigenvectors of $M$ corresponding to the eigenvalues $\lambda_1(M)$ and $\lambda_2(M)$, respectively.
Then
\[
 \frac{\Pr(\lambda_2(W)\ge x)}{\Pr(\lambda_1(W)\ge x)}
 \le \frac{\Pr\bigl(\sigma\chi^2(n;\lambda_2(M)/\sigma)\ge x\bigr)}{\Pr\bigl(\sigma\chi^2(n;\lambda_1(M)/\sigma)\ge x\bigr)},
\]
which goes to zero as $x$ goes to infinity.
\hfill\qedsymbol


\section{Central case with a scalar covariance: Selberg type integral and Laguerre polynomials}
\label{sec:central}

We assume that $M=0$ (central) and 
$\Sigma$ in (\ref{parameters}) is a scalar matrix, and  we study this case by special functions.
Under these assumptions, 
we show that the expectation of the Euler characteristic 
can be expressed in terms of a Selberg type integral, which is equal to
a Laguerre polynomial
in view of the works by Aomoto~\cite{Aomoto} and Kaneko~\cite{Kaneko}.

\begin{theorem}
\label{th:e_is_selberg_integral}
Let
$$ M_x = \{ hg^{\top}\,|\, g^{\top} A h \geq x, h \in S^{n-1}, g \in S^{m-1} \}, \quad m \leq n. $$
Assume that the distribution of $m \times n$ random matrices $A$ 
is the Gaussian distribution with mean $0$ and covariance $I_m/s$.
In other words, we have
$$ p(A) \sim \exp\left( -\frac{1}{2} \tr\,(sA^{\top} A) \right). $$
Then we have
\[
{\rm E}[\chi(M_x(s))] = 
\prod_{i=1}^5 c_i \int_x^{+\infty} \exp\left( - \frac{s}{2} \sigma^2\right)
 \sigma^{n-m} {}_1 F_1(-(m-1),1+n-m;s \sigma^2) d \sigma,
\]
where $c_1, c_2, c_3, c_4, c_5$ are given by
(\ref{eq:c1}), (\ref{eq:c2}), (\ref{eq:c3}), (\ref{eq:c4}), (\ref{eq:c5}), respectively.
\end{theorem}

\textbf{Proof.}
For $g \in S^{m - 1}, h \in S^{n-1}$, set
$$
 {\widetilde G} =
\Bigl( g \ \big| \ G \Bigr)
 \in O(m), \quad
\mbox{$g$ is a column vector,}
$$
$$
 {\widetilde H} =
\Bigl( h \ \big| \ H \Bigr)
 \in O(n), \quad
\mbox{$h$ is a column vector.}
$$
Then the $m \times n$ matrix $A$ can be written as
$$ A = {\widetilde G}
\left(
 \begin{array}{c|c}
  \sigma & 0 \\ \hline
  0 & B
 \end{array}
\right)
{\widetilde H}^{\top}. 
$$ 
We denote by ${\widetilde B}$ the middle matrix in the above expression.

Set $\etr(X) = \exp(\tr(X))$ and $S=\Sigma^{-1}$.
We consider the central case $M=0$ in (\ref{eq:chi-for-normal}).
Since $\tr(PQ)=\tr(QP)$ and ${\widetilde H}^{\top} {\widetilde H} =E$, 
we have
\[
 \etr\left( -\frac{1}{2} A^{\top} S A \right)
 = \etr\left( -\frac{1}{2} {\widetilde H} {\tilde B}^{\top} {\widetilde G}^{\top} S {\widetilde G} {\widetilde B} {\widetilde H}^{\top} \right) 
 = \etr\left( -\frac{1}{2} S {\widetilde G} {\widetilde B} {\widetilde H}^{\top} {\widetilde H} {\widetilde B}^{\top} {\widetilde G}^{\top} \right) 
 = \etr\left( -\frac{1}{2} S {\widetilde G} ({\widetilde B} {\widetilde B}^{\top}) {\widetilde G}^{\top} \right).
\] 
It follows from Theorem \ref{th:expectation_of_euler_characteristic} with $p(A)$ being the normal distribution that 
\begin{align*}
{\rm E}[\chi(M_x)]
 =& c_1(m,n,s) \int_x^\infty \sigma^{n-m} d \sigma
 \int_{\mathbb{R}^{(m-1)(n-1)}} dB
 \int_{S^{m-1}} G^{\top} dg 
 \int_{S^{n-1}} H^{\top} dh \\
&\times
\det\left(\sigma^2 I_{m-1} - B B^{\top}\right) 
 \etr\left( -\frac{1}{2} S {\widetilde G} ({\widetilde B} {\widetilde B}^{\top}) {\widetilde G}^{\top} \right),
\end{align*}
where
\begin{equation}
\label{eq:c1}
  c_1(m,n,s) = \frac{1}{2} \cdot
  \frac{1}{(2 \pi)^{nm/2} \det(S^{-1})^{n/2}},
 \quad
 S = s I_m.
\end{equation}
We denote by $G_i$ the $i$-th column vector of $G$ and by $dg$
the column vector of the differential form $dg_i$.
Define
$$ G^{\top} dg = \wedge_{i=1}^{m-1} G_i^{\top} dg. $$
It is an invariant measure for rotations on $S^{m-1}$
\cite[Theorem 4.2]{James1954}.
We may define $H^{\top} dh$ analogously.

Moreover, since $S = s I_m$, we have 
$$
 \etr\left( -\frac{1}{2} S {\widetilde G} ({\widetilde B} {\widetilde B}^{\top}) {\widetilde G}^{\top} \right) 
 = \etr\left( -\frac{s}{2} {\widetilde B} {\widetilde B}^{\top} \right).
$$
Since there is no $G, H$ involved on the right side of the above identity,
we can separate the following integral:
\begin{equation}
\label{eq:c2}
c_2(m,n)=
\int_{S^{m-1}} G^{\top} dg 
\int_{S^{n-1}} H^{\top} dh 
= \frac{2 \pi^{m/2}}{\Gamma(m/2)}
  \frac{2 \pi^{n/2}}{\Gamma(n/2)}.
\end{equation}
Therefore, we only need to evaluate the integral
\[
 \int_{\mathbb{R}^{(m-1)(n-1)}} dB \,
 \det\left(\sigma^2 I_{m-1} - B B^{\top}\right) 
\etr\left(-\frac{s}{2} {\widetilde B} {\widetilde B}^{\top}\right).
\]
We denote the integral above by $\qq(s; \sigma)$.
In terms of $\qq(s;\sigma)$, we have
$$
 {\rm E}[\chi(M_x)] = c_1(S) c_2(m) \int_x^\infty \sigma^{n-m} \qq(s;\sigma) d \sigma.
$$

We make the singular value decomposition of the matrix $B$ as
$B=PLQ^{\top}$, where the matrices $P \in O(m-1)$, 
 $Q \in V_{m-1}(\mathbb{R}^{n-1})$ (Stiefel manifold), 
$L=\diag(\ell_1, \ldots, \ell_{m-1})$
(see, \eg, \cite{James1954} and~\cite[(3.1)]{takemura-kuriki-1999}).
It follows from \cite[(3.1)]{takemura-kuriki-1999}
that
$$ dB = \prod_{i=1}^{m-1} \ell_i^{(n-1)-(m-1)}
   \prod_{1 \leq i < j \leq m-1} (\ell_i^2-\ell_j^2) 
   \left( \prod_{i=1}^{m-1} d\ell_i \right) \wedge
   \omega_1 \wedge \omega_2,
$$
$$ \omega_1 = 
   \wedge_{1 \leq i \leq m-1, i < j \leq m-1} P_j^{\top} dP_i,
\qquad
 \omega_2 = \wedge_{1 \leq i \leq m-1, i < j \leq n-1} Q_j^{\top} dQ_i
$$
the volume element of the Stiefel manifold,
when $\ell_1 \geq \ell_2 \geq \cdots \geq \ell_{m-1} \geq 0$.
Here, $P_i$ is the $i$-th column vector of $P$.
Since
$$
 \det\left( \sigma^2 I_{m-1} - PLQ^{\top} Q L^{\top} P^{\top}\right)
 = \det\left(P(\sigma^2 I_{m-1} - L L^{\top})P^{\top}\right) 
 = \det\left(\sigma^2 I_{m-1} - L L^{\top}\right),
$$
and
\[
\etr\left(-\frac{s}{2} {\tilde B} {\tilde B}^{\top}\right)
= \exp\left(-\frac{s}{2} \sigma^2\right)
  \etr\left(-\frac{s}{2} B B^{\top}\right) 
= \exp\left(-\frac{s}{2} \sigma^2\right)
  \etr\left(-\frac{s}{2} PL Q^{\top} Q L^{\top} P^{\top}\right) 
= \exp\left(-\frac{s}{2} \sigma^2\right)
  \exp\left(-\frac{s}{2} L L^{\top} \right),
\]
 we have
\begin{equation}
\label{eq:sel1}
\qq(s;\sigma) = c_3'(m,n;\sigma) \int_{L \in \mathbb{R}^{m-1}}
   \prod_{i=1}^{m-1} |\ell_i|^{n-m}
   \prod_{1 \leq i < j \leq m-1} |\ell_i^2-\ell_j^2| 
   \prod_{i=1}^{m-1} (\sigma^2-\ell_i^2)
    \exp\left(-\frac{s}{2} \sum \ell_i^2 \right) 
   \prod_{i=1}^{m-1} d\ell_i,
\end{equation}
where $c_3'(m,n;\sigma)=c_3(m,n,;\sigma) \exp\left(-\frac{s}{2} \sigma^2 \right)$,
\begin{align}
c_3(m,n;\sigma)
=& \frac{1}{(m-1)!2^{m-1} 2^{m-1}}
   \int_{O(m-1)} \omega_1 \int_{V_{m-1}(\mathbb{R}^{n-1})} \omega_2 \nonumber \\
=& \frac{1}{(m-1)! 2^{m-1} 2^{m-1}}
   2^{m-1} \prod_{k=1}^{m-1} \frac{\pi ^{k/2}}{\Gamma(k/2)} 
  \times
   \frac{2^{m-1} \pi^{(m-1)(n-1-(m-2)/2)/2}}
        {\prod_{i=1}^{m-1} \Gamma((n-1)/2-(i-1)/2)}.
\label{eq:c3}
\end{align}
In~\eqref{eq:c3}, there is a constant $(m-1)! 2^{m-1} 2^{m-1}$ involved in the denominator because in this case~$(m-1)! 2^{m-1}$ copies of 
the domain $\ell_1 \geq \ell_2 \geq \ldots \geq \ell_{m-1} \geq 0$
cover $\mathbb{R}^{m-1}$, and the correspondence between the coordinates of $B$ and those of its singular value decomposition is $1 / 2^{m-1}$
because we have the choice of signs of the eigenvector $P_i$.
For the volumes of $O(m-1)$ and $V_{m-1}(\mathbb{R}^{n-1})$, see, \eg, \cite[Proposition 2.23, Theorem 2.24]{woods-zhang-2015}.

In (\ref{eq:sel1}),
we make a change of variables by
$\ell_i' = \ell_i^2$.
Then we have
$d \ell_i' = 2 \ell_i d \ell_i$,
and 
$$ d \ell_i = \frac{1}{2 \sqrt{\ell_i'}} d \ell_i'. $$
Furthermore, we have
\[
\qq(s;\sigma) = c'_3(m,n;\sigma) \nonumber \int_{L' \in \mathbb{R}_{\geq 0}^{m-1}}
  \prod_{i=1}^{m-1} {\ell_i'}^{-1/2+(n-m)/2} 
  \prod_{1 \leq i < j \leq m-1} |\ell_i'-\ell_j'| 
  \prod_{i=1}^{m-1} (\sigma^2-\ell_i') 
   \exp\left(-\frac{s}{2} \sum \ell_i' \right) 
  \prod_{i=1}^{m-1} d\ell_i'.  
\]
Set $\ell_i' = \frac{2}{s} \ell_i''$ and factor out $s >0$.
Then it follows from $d \ell_i' = \frac{2}{s} d\ell_i''$ that 
$$
 \qq(s;\sigma)= c'_3(m,n;\sigma) c_4(m,n,s) {\tilde \qq}(s;\sigma),
$$
where 
\begin{equation} 
{\tilde \qq}(s;\sigma)
 = \int_{L'' \in \mathbb{R}_{\geq 0}^{m-1}}
 \prod_{i=1}^{m-1} {\ell_i''}^{-1/2+(n-m)/2} 
   \prod_{1 \leq i < j \leq m-1} |\ell_i''-\ell_j''| 
   \prod_{i=1}^{m-1} \left(\frac{\sigma^2 s}{2}-\ell_i''\right)
    \exp\left(-\sum_{i=1}^{m-1} \ell_i'' \right) 
   \prod_{i=1}^{m-1} d\ell_i'', 
\label{eq:sel2b}
\end{equation}
and
\begin{equation}
\label{eq:c4}
c_4(m,n,s)
 = (s/2)^{(m-1)/2} (s/2)^{(m-n)(m-1)/2} (s/2)^{-\frac{1}{2}(m-1)(m-2)} (s/2)^{-(m-1)} (s/2)^{-(m-1)}
 = (s/2)^{-\frac{1}{2}(m^2-1)-\frac{1}{2}(n-m)(m-1)}.
\end{equation}
The integral~\eqref{eq:sel2b} can be expressed as a polynomial in $\sigma$.
Let us derive differential equations for this integral and express it
in terms of a special polynomial.
We utilize the result by Aomoto \cite{Aomoto} and its generalization
by Kaneko~\cite{Kaneko}.
In \cite{Kaneko}, a system of differential equations,
special values, and an expansion in terms of Jack polynomials
were given 
for the integral 
\begin{equation}
\int_{[0,1]^{m-1}}
\prod_{1 \leq i \leq m-1, 1\leq k \leq r} (\ell_i-\sigma_k)^\mu
D(\ell_1, \ldots, \ell_{m-1}) d\ell_1 \cdots d\ell_{m-1},
\label{eq:kaneko93}
\end{equation}
$$ D =
 \prod_{i=1}^{m-1} \ell_i^{\lambda_1} (1-\ell_i)^{\lambda_2}
 \prod_{1 \leq i < j \leq m-1} | \ell_i - \ell_j |^\lambda,
$$
when $\mu=1$ or $\mu = -\lambda/2$.
Let us make the coordinate change
$ \ell_i = y_i/N $, $\lambda_2=N$, $\sigma_i = \tau_i/N$.
Then we have
$d \ell_i = dy_i/N$,
$(1-\ell_i)^\lambda = (1-y_i/N)^N$,
$$ (1-y_i/N)^N \rightarrow \exp(-y_i), \quad N \rightarrow \infty. $$
The integral (\ref{eq:kaneko93}) becomes 
$$
c_N \int_{[0,N]^{m-1}}
\prod_{1 \leq i \leq m-1, 1\leq k \leq r} (y_i-\tau_k)^\mu
D(y_1, \ldots, y_{m-1}) dy_1 \cdots dy_{m-1},
$$
$$
 D =
 \prod_{i=1}^{m-1} y_i^{\lambda_1} (1-y_i/N)^N
 \prod_{1 \leq i < j \leq m-1} | y_i - y_j |^\lambda,
 c_N = N^{-r(m-1)-(m-1)-\lambda_1 (m-1) - \lambda (m-1)(m-2)/2}. 
$$
When $N \rightarrow \infty$,
the above integral divided by $c_N$ converges to
\begin{equation}
\label{EQ:converge}
\int_{\mathbb{R}_{\geq 0}^{m-1}}
\prod_{1 \leq i \leq m-1, 1\leq k \leq r} (y_i-\tau_k)^\mu
D(y_1, \ldots, y_{m-1}) dy_1 \cdots dy_{m-1},
\end{equation}
$$
 D =
 \prod_{i=1}^{m-1} y_i^{\lambda_1} \exp\left(-\sum_{i=1}^{m-1} y_i\right)
 \prod_{1 \leq i < j \leq m-1} | y_i - y_j |^\lambda.
$$
Let us apply this limiting procedure to derive a differential
equation for the above integral.
When $r= \mu=1$,
the differential equation for the integral (\ref{eq:kaneko93}) is
\begin{equation}
\label{EQ:gauss}
 \sigma(1-\sigma) \pd{\sigma}^2+(c-(a+b+1)\sigma) \pd{\sigma} - ab,
\end{equation}
where
$ a= -(m-1)$, 
$b=\frac{2}{\lambda}(\lambda_1+\lambda_2+2)+(m-1)+1$,
$c = \frac{2}{\lambda}(\lambda_1+1)$.
This is the Gauss hypergeometric equation.
Set
$\lambda_2=N$, $\sigma = \frac{z}{N}$.
Then we can find the limit of this equation when $N \rightarrow \infty$.
In fact, it can be performed as follows.
Set $\theta_z = z \pd{z}$. Note that~\eqref{EQ:gauss} is invariant by the scalar multiplication of $z$.
Then the limit of
$$
 \theta_z (\theta_z + \frac{2}{\lambda}(\lambda_1+1)-1)
  - \frac{z}{N} (\theta_z-(m-1))(\theta_z+\frac{2}{\lambda}(N+\lambda_1+2)+(m-1)+1),
$$
when $N \rightarrow \infty$, equals

$$
 \theta_z (\theta_z + \frac{2}{\lambda}(\lambda_1+1)-1)
  - \frac{2}{\lambda} z (\theta_z-(m-1)).
$$
In particular, when $\lambda=1$ and $\lambda_1=-1/2+(n-m)/2$, 
 we have
$$ \theta_z (\theta_z + n - m)  - 2 z (\theta_z-(m-1)). $$
A polynomial solution of the above equation can be written as 
a constant multiple of the confluent hypergeometric
polynomial ${}_1F_1(-(m-1),1+n-m; 2z)$.
Therefore, it follows from~\eqref{eq:sel2b},~\eqref{EQ:converge} and the above argument that
\begin{align*}
\qq(s;\sigma)
=& c_3'(m,n;\sigma) c_4(m,n,s) c_5(m,n) \cdot {}_1F_1(-(m-1),1+n-m; \sigma^2s) \\
=& c_3'(m,n;\sigma) c_4(m,n,s) c_5(m,n) \left(
1 
+ \frac{-(m-1)}{1 (n-m+1)} (\sigma^2s)
+ \frac{(m-1) (m-2) }{2! (n-m+1)(n-m+2)} (\sigma^2s)^2 \right. \\
& \quad \left.
+ \frac{-(m-1) (m-2)(m-3) }{3!(n-m+1)(n-m+2)(n-m+3)} (\sigma^2 s)^3
+ \cdots \right. \nonumber \\
& \quad \left. + \frac{(-1)^{m-1} (m-1)!}{(m-1)!(n-m+1) \cdots (n-m+1+m-1)} (\sigma^2 s)^{m-1} \right),
\end{align*}
where
\begin{equation}
\label{eq:c5}
 c_5(m,n) = \mbox{(the expression (\ref{eq:sel2b}))}|_{\sigma=0}
 = \prod_{i=1}^{m-1} \frac{\Gamma\left(1+\frac{i}{2}\right)\Gamma\left(\frac{3}{2}+\frac{n-m}{2}+\frac{i-1}{2}\right)}
 {\Gamma\left(\frac{3}{2}\right)} (-1)^{m-1}
\end{equation}
by taking a limit of the Selberg integral formula
\cite{selberg}
with an analogous method as was used when deriving (\ref{EQ:converge}).
\qedsymbol

\begin{figure}[htb]
\begin{center}
\includegraphics[width=6cm]{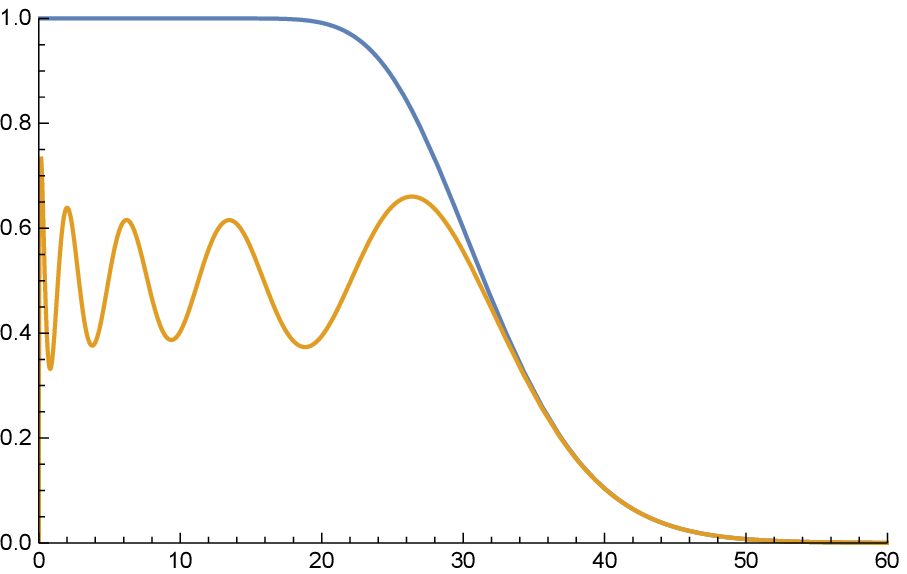} \quad
\includegraphics[width=6cm]{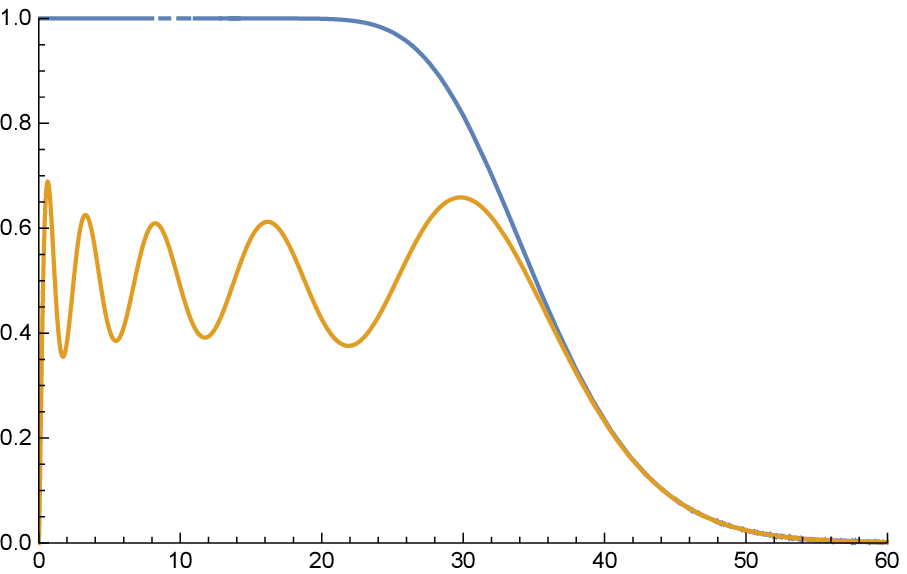}
\end{center}
\caption{Approximation versus exact values of ${\rm E}[\chi(M_x)]$: $m=n=10$, $s=1$ (left) and $m=10, n=12$, $s=1$ (right).}
\label{fig:approx:m10}
\end{figure}
\begin{figure}[htb]
\begin{center}
\includegraphics[width=6cm]{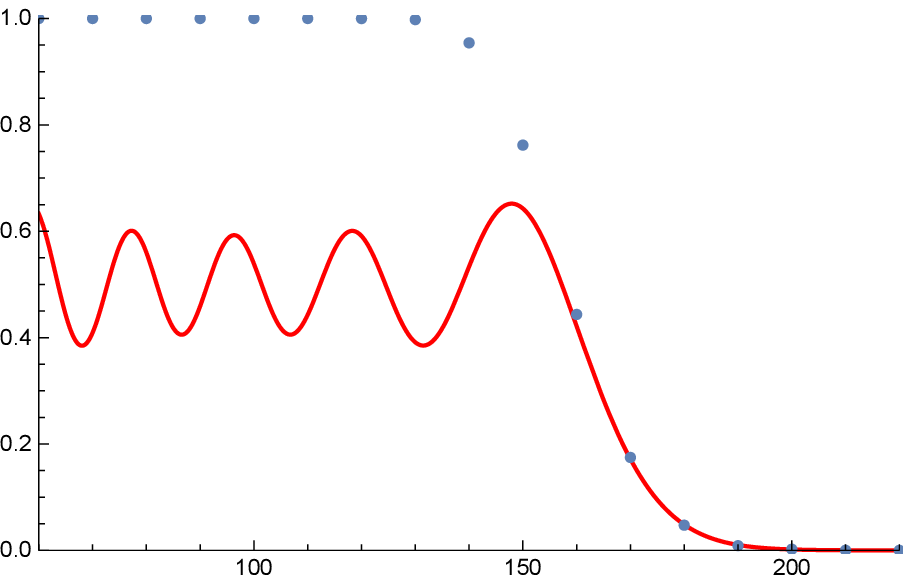} \quad
\includegraphics[width=6cm]{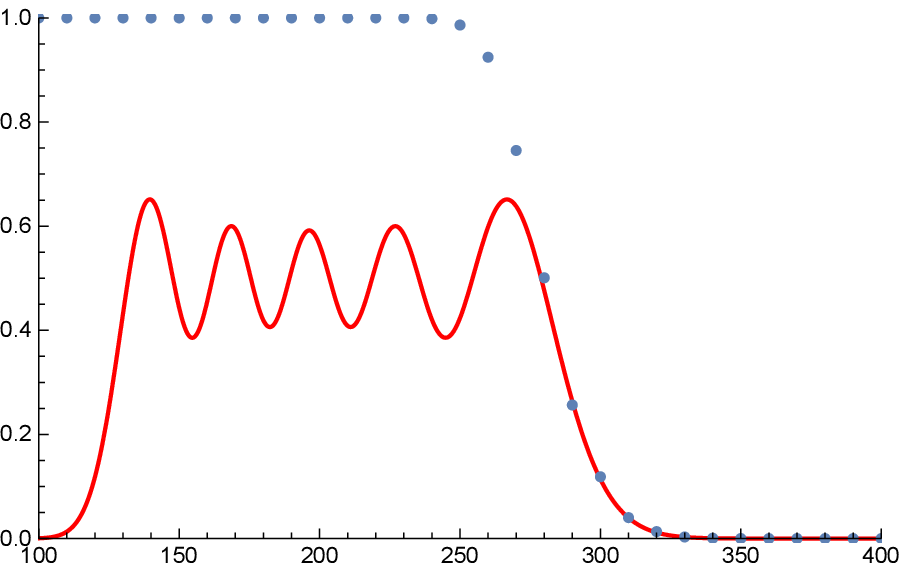}
\end{center}
\caption{Approximation (curves) vs simulation (dots) of ${\rm E}[\chi(M_x)]$: $m=10, n=100$, $s=1$ (left) and $m=10, n=200$, $s=1$ (right).}
\label{fig:approx:m10n100}
\end{figure}

Let us make a numerical evaluation by utilizing Theorem 
\ref{th:e_is_selberg_integral} when $m=n=3$.
If $m=n=3$, we have
$$
 c_1 c_2 c_3 c_4 c_5 = 2 \sqrt{2/\pi} \sqrt{s}.
$$
Since
\[
u(s,k,x) = \int_x^{+\infty} \exp(-\sigma^2 s/2) \sigma^{2k} d\sigma 
= \Gamma(k+1/2) \left( \frac{2}{s} \right)^{k+1/2} \frac{1}{2}
\int_{x^2}^{+\infty}
\frac{ y^{k+1/2-1} \exp(-y/(2/s)) dy}{\Gamma(k+2) (2/s)^{k+1/2}}, 
\]
where the last integral is equal to the upper tail probability
of the Gamma distribution with scale parameter $2/s$ and shape parameter $k+1/2$.
It follows from Theorem \ref{th:e_is_selberg_integral} that 
the expectation
is equal to 
\begin{equation}
\label{EQ:eulercentral}
 {\rm E}[\chi(M_x)] = 2 \sqrt{2/\pi} \sqrt{s} \left( u(s,0,x)-2 s u(s,1,x) + \frac{s^2}{2} u(s,2,x) \right).
\end{equation}
An \textsf{R} code for evaluating ${\rm E}[\chi(M_x)]$ in this case is as follows:
{\footnotesize
\begin{verbatim}
ug2<-function(s,k,x) {
  return(pgamma(x^2, scale=2/s, shape=k+1/2, lower = FALSE)*
        gamma(k+1/2)*(2/s)^(k+1/2)/2);
}
ec3<-function(x,s) {
 cc<- 2*(2/pi)^(1/2)*s^(1/2);
 c5<-1;
 return(cc*c5*
   (ug2(s,0,x)-2*s*ug2(s,1,x)+(1/2)*s^2*ug2(s,2,x)));
}
## Draw a graph
curve(ec3(x,1),from=1,to=10)
\end{verbatim}
}
When $s=1$, some values are given in Table~\ref{TABLE:5}:
\begin{table}[h]
\caption{Evaluation of ${\rm E}[\chi(M_x)]$ by~\eqref{EQ:eulercentral} versus Monte Carlo simulation}
\begin{center}
\begin{tabular}{ccc}
\hline
$x$ & ${\rm E}[\chi(M_x)]$ & simulation (with 100000 tries) \\ \hline
3   & 0.215428520          & 0.217072 \\
4   & 0.016122970          & 0.016195 \\
5   & 0.000357368          & 0.000386 \\
\hline
\end{tabular}
\end{center}
\label{TABLE:5}
\end{table}
We present two graphs in Fig.~\ref{fig:approx:m10} to compare our approximate formula with the exact values by the Pfaffian of a matrix
(see, \eg,~\cite{KC1971, Chiani2016}).
The matrix sizes are $m=n=10$ and $m=10, n=12$,
and $s=1$.
The horizontal axis is $x^2$.
Note that our approximation formula is expressed as a finite sum of 
$m$ terms of incomplete Gamma functions which can be evaluated
faster than the Pfaffian of an $m \times m$ matrix when $m$ becomes larger.
The approximation error was evaluated by \cite{kuriki-takemura-2001, kuriki-takemura-2008} as
\[
 \Delta(x)=\mathrm{E}[\chi(M_x)]-\Pr(\lambda_1(W)\ge x^2) \sim
 -\frac{1}{\Gamma(m-1)\Gamma(n-1)}x^{2(m+m-5)}e^{-x^2}, \quad x\to\infty,
\]
which is exponentially smaller than
\[
 \mathrm{E}[\chi(M_x)]
 \sim \Pr(\lambda_1(W)\ge x^2)
 \sim \frac{\sqrt{\pi}}{2^{(m+n-5)/2}\Gamma(m/2)\Gamma(n/2)} x^{m+n-3} e^{-x^2/2}, \quad x\to\infty.
\]
This explains the very accurate tail behaviors in Figs.~\ref{fig:approx:m10} and \ref{fig:approx:m10n100}.
Note that $\Delta(x)$ is always negative because (\ref{altsum}) is always less than $\Pr(\lambda_1(W)\ge x)$. 

The two graphs in Fig.~\ref{fig:approx:m10n100} are to compare our approximate formula with values by a simulation of $10000$ tries.
The matrix size is $m=10, n=100$ and $m=10, n=200$, respectively,
and $s=1$.
The horizontal axis is $x^2$. 

\section*{Acknowledgments}

The authors would like to thank the Editor-in-Chief, Associate Editors and two anonymous referees who kindly reviewed 
the earlier versions of this paper and provided valuable comments and suggestions. 
Besides, we also deeply thank Christoph Koutschan, who is the author
of the package \texttt{HolonomicFunctions} used in this study,
for his help and encouragement. 
This research is partially supported by the Austrian Science Fund (FWF): P29467-N32, JSPS KAKENHI Grant Number 16H02792, 
 the UTD start-up grant: P-1-03246, the Natural Science Foundation of USA grants: CCF-1815108 and CCF-1708884, and JST CREST Grant Number JP19209317.

\section*{References}

\bibliographystyle{myjmva}  

\begin{thebibliography}{38}
\expandafter\ifx\csname natexlab\endcsname\relax\def\natexlab#1{#1}\fi
\providecommand{\bibinfo}[2]{#2}
\ifx\xfnm\relax \def\xfnm[#1]{\unskip,\space#1}\fi
\bibitem[{Adler(1981)}]{adler-1981}
\bibinfo{author}{R.~J. Adler}, \bibinfo{title}{The geometry of random fields},
  \bibinfo{publisher}{John Wiley \& Sons, Ltd., Chichester},
  \bibinfo{year}{1981}. \bibinfo{note}{Wiley Series in Probability and
  Mathematical Statistics}.
\bibitem[{Adler and Taylor(2007)}]{Adler2007}
\bibinfo{author}{R.~J. Adler}, \bibinfo{author}{J.~E. Taylor},
  \bibinfo{title}{Random fields and geometry}, \bibinfo{publisher}{Springer},
  \bibinfo{address}{New York}, \bibinfo{year}{2007}.
\bibitem[{Aomoto(1987)}]{Aomoto}
\bibinfo{author}{K.~Aomoto}, \bibinfo{title}{Jacobi polynomials associated with
  selberg integrals}, \bibinfo{journal}{SIAM Journal on Mathematical Analysis}
  \bibinfo{volume}{18} (\bibinfo{year}{1987}) \bibinfo{pages}{545--549}.
\bibitem[{Beran(1975)}]{beran-1975}
\bibinfo{author}{R.~Beran}, \bibinfo{title}{Tail probabilities of noncentral
  quadratic forms}, \bibinfo{journal}{The Annals of Statistics}
  \bibinfo{volume}{3} (\bibinfo{year}{1975}) \bibinfo{pages}{969--974}.
\bibitem[{Chen et~al.(2019)Chen, Kauers, Li and Zhang}]{Yi2017}
\bibinfo{author}{S.~Chen}, \bibinfo{author}{M.~Kauers},
  \bibinfo{author}{Z.~Li}, \bibinfo{author}{Y.~Zhang}, \bibinfo{title}{Apparent
  singularities of {D}-finite systems}, \bibinfo{journal}{Journal of Symbolic
  Computation} \bibinfo{volume}{95} (\bibinfo{year}{2019})
  \bibinfo{pages}{217--237}.
\bibitem[{Chiani(2016)}]{Chiani2016}
\bibinfo{author}{M.~Chiani}, \bibinfo{title}{Distribution of the largest root
  of a matrix for roy's test in multivariate analysis of variance},
  \bibinfo{journal}{Journal of Multivariate Analysis} \bibinfo{volume}{143}
  (\bibinfo{year}{2016}) \bibinfo{pages}{467--471}.
\bibitem[{Chikuse(1992)}]{Chikuse1992}
\bibinfo{author}{Y.~Chikuse}, \bibinfo{title}{Properties of {H}ermite and
  {L}aguerre polynomials in matrix argument and their applications},
  \bibinfo{journal}{Linear Algebra and its Applications} \bibinfo{volume}{176}
  (\bibinfo{year}{1992}) \bibinfo{pages}{237--260}.
\bibitem[{Chyzak(2000)}]{Chyzak2000}
\bibinfo{author}{F.~Chyzak}, \bibinfo{title}{An extension of {Z}eilberger's
  fast algorithm to general holonomic functions}, \bibinfo{journal}{Discrete
  Mathematics} \bibinfo{volume}{217 (1-3)} (\bibinfo{year}{2000})
  \bibinfo{pages}{115--134}.
\bibitem[{Coxe et~al.(2015)Coxe, Little and O'Shea}]{Cox2015}
\bibinfo{author}{D.~A. Coxe}, \bibinfo{author}{J.~Little},
  \bibinfo{author}{D.~O'Shea}, \bibinfo{title}{Ideals, Varieties, and
  Algorithms}, \bibinfo{publisher}{Springer}, \bibinfo{address}{New York},
  \bibinfo{edition}{4th} edition, \bibinfo{year}{2015}.
\bibitem[{Danufane et~al.(2017)Danufane, Siriteanu, Ohara and
  Takayama}]{Danufane2017}
\bibinfo{author}{F.~H. Danufane}, \bibinfo{author}{C.~Siriteanu},
  \bibinfo{author}{K.~Ohara}, \bibinfo{author}{N.~Takayama},
  \bibinfo{title}{Holonomic gradient method-based cdf evaluation for the
  largest eigenvalue of a complex noncentral {W}ishart matrix},
  \bibinfo{year}{2017}. \bibinfo{note}{ArXiv:1707.02564}.
\bibitem[{Davis(1979)}]{Davis1979}
\bibinfo{author}{A.~W. Davis}, \bibinfo{title}{Invariant polynomials with two
  matrix arguments extending the zonal polynomials: Applications to
  multivariate distribution theory}, \bibinfo{journal}{Annals of the Institute
  of Statistical Mathematics} \bibinfo{volume}{31} (\bibinfo{year}{1979})
  \bibinfo{pages}{465--485}.
\bibitem[{Davis(1980)}]{Davis1980}
\bibinfo{author}{A.~W. Davis}, \bibinfo{title}{Invariant polynomials with two
  matrix arguments: extending the zonal polynomials}, in:
  \bibinfo{editor}{P.~R. Krishnaiah} (Ed.), \bibinfo{booktitle}{Multivariate
  Analysis V}, \bibinfo{publisher}{North-Holland Publishing Company},
  \bibinfo{year}{1980}, pp. \bibinfo{pages}{287--299}.
\bibitem[{D\'{i}az-Garc\'{i}a and Guti\'{e}rrez-J\'{a}imez(2011)}]{Garcia2011}
\bibinfo{author}{J.~A. D\'{i}az-Garc\'{i}a},
  \bibinfo{author}{R.~Guti\'{e}rrez-J\'{a}imez}, \bibinfo{title}{On {W}ishart
  distribution: Some extensions}, \bibinfo{journal}{Linear Algebra and its
  Applications} \bibinfo{volume}{435} (\bibinfo{year}{2011})
  \bibinfo{pages}{1296--1310}.
\bibitem[{Hashiguchi et~al.(2013)Hashiguchi, Numata, Takayama and
  Takemura}]{Hashiguchi2013}
\bibinfo{author}{H.~Hashiguchi}, \bibinfo{author}{Y.~Numata},
  \bibinfo{author}{N.~Takayama}, \bibinfo{author}{A.~Takemura},
  \bibinfo{title}{Holonomic gradient method for the distribution function of
  the largest root of a wishart matrix}, \bibinfo{journal}{Journal of
  Multivariate Analysis} \bibinfo{volume}{117} (\bibinfo{year}{2013})
  \bibinfo{pages}{296--312}.
\bibitem[{Hayakawa(1969)}]{Hayakawa1969}
\bibinfo{author}{T.~Hayakawa}, \bibinfo{title}{On the distribution of the
  latent roots of a positive definite random symmetric matrix i},
  \bibinfo{journal}{Annals of the Institute of Statistical Mathematics}
  \bibinfo{volume}{21} (\bibinfo{year}{1969}) \bibinfo{pages}{1--21}.
\bibitem[{Hibi and {\rm et~al.}(2013)}]{Hibi2013}
\bibinfo{author}{T.~Hibi}, \bibinfo{author}{{\rm et~al.}},
  \bibinfo{title}{Gr\"{o}bner Bases: Statistics and software systems},
  \bibinfo{publisher}{Springer}, \bibinfo{address}{New York},
  \bibinfo{year}{2013}.
\bibitem[{Hillebrand and Schmale(2001)}]{Hillebrand2001}
\bibinfo{author}{A.~Hillebrand}, \bibinfo{author}{W.~Schmale},
  \bibinfo{title}{Towards a effective version of a theorem of {S}tafford},
  \bibinfo{journal}{Journal of Symbolic Computation} \bibinfo{volume}{32}
  (\bibinfo{year}{2001}) \bibinfo{pages}{699--716}.
\bibitem[{James(1954)}]{James1954}
\bibinfo{author}{A.~T. James}, \bibinfo{title}{Normal multivariate analysis and
  the orthogonal group}, \bibinfo{journal}{The Annals of Mathematical
  Statistics} \bibinfo{volume}{25} (\bibinfo{year}{1954})
  \bibinfo{pages}{40--75}.
\bibitem[{James(1955)}]{James1955}
\bibinfo{author}{A.~T. James}, \bibinfo{title}{The non-central {W}ishart
  distribution}, \bibinfo{journal}{Proceedings of the Royal Society of London}
  \bibinfo{volume}{229} (\bibinfo{year}{1955}) \bibinfo{pages}{364--366}.
\bibitem[{Kaneko(1993)}]{Kaneko}
\bibinfo{author}{J.~Kaneko}, \bibinfo{title}{Selberg integrals and
  hypergeometric functions associated with {J}ack polynomials},
  \bibinfo{journal}{SIAM Journal on Mathematical Analysis} \bibinfo{volume}{24}
  (\bibinfo{year}{1993}) \bibinfo{pages}{1086--1110}.
\bibitem[{Kang and Alouini(2003)}]{Kang2003}
\bibinfo{author}{M.~Kang}, \bibinfo{author}{M.~S. Alouini},
  \bibinfo{title}{Largest eigenvalue of complex wishart matrices and
  performance analysis of {MIMO} {MRC} systems}, \bibinfo{journal}{IEEE Journal
  on Selected Areas in Communications} \bibinfo{volume}{21}
  (\bibinfo{year}{2003}) \bibinfo{pages}{418--426}.
\bibitem[{Kauers et~al.(2009)Kauers, Koutschan and Zeilberger}]{KKZ2009}
\bibinfo{author}{M.~Kauers}, \bibinfo{author}{C.~Koutschan},
  \bibinfo{author}{D.~Zeilberger}, \bibinfo{title}{Proof of {I}ra {G}essel's
  lattice path conjecture}, \bibinfo{journal}{Proceedings of the National
  Academy of Sciences} \bibinfo{volume}{106 (28)} (\bibinfo{year}{2009})
  \bibinfo{pages}{1150211505}.
\bibitem[{Koutschan(2009)}]{Christoph2009}
\bibinfo{author}{C.~Koutschan}, \bibinfo{title}{Advanced applications of the
  holonomic systems approach}, Ph.D. thesis, Johannes Kepler University Linz,
  \bibinfo{year}{2009}.
\bibitem[{Koutschan(2010)}]{Christoph2010}
\bibinfo{author}{C.~Koutschan}, \bibinfo{title}{HolonomicFunctions user's
  guide}, \bibinfo{type}{Technical Report}, Johannes Kepler University Linz,
  \bibinfo{year}{2010}.
  \bibinfo{note}{\href{http://www.risc.jku.at/publications/download/risc_3934/hf.pdf}{http:/$\!$/www.risc.jku.at/publications/download/risc\_3934/hf.pdf}}.
\bibitem[{Koutschan et~al.(2011)Koutschan, Kauers and Zeilberger}]{KKZ2011}
\bibinfo{author}{C.~Koutschan}, \bibinfo{author}{M.~Kauers},
  \bibinfo{author}{D.~Zeilberger}, \bibinfo{title}{Proof of {G}eorge
  {A}ndrews's and {D}avid {R}obbins's $q$-{TSPP} conjecture},
  \bibinfo{journal}{Proceedings of the National Academy of Sciences}
  \bibinfo{volume}{108 (6)} (\bibinfo{year}{2011}) \bibinfo{pages}{21962199}.
\bibitem[{Krishnaiah and Chang(1971)}]{KC1971}
\bibinfo{author}{P.~R. Krishnaiah}, \bibinfo{author}{T.~C. Chang},
  \bibinfo{title}{On the exact distributions of the extreme roots of the
  wishart and manova matrices}, \bibinfo{journal}{Journal of Multivariate
  Analysis} \bibinfo{volume}{1} (\bibinfo{year}{1971})
  \bibinfo{pages}{108--117}.
\bibitem[{Kuriki and Takemura(2001)}]{kuriki-takemura-2001}
\bibinfo{author}{S.~Kuriki}, \bibinfo{author}{A.~Takemura},
  \bibinfo{title}{Tail probabilities of the maxima of multilinear forms and
  their applications}, \bibinfo{journal}{The Annals of Statistics}
  \bibinfo{volume}{29} (\bibinfo{year}{2001}) \bibinfo{pages}{328--371}.
\bibitem[{Kuriki and Takemura(2008)}]{kuriki-takemura-2008}
\bibinfo{author}{S.~Kuriki}, \bibinfo{author}{A.~Takemura},
  \bibinfo{title}{Euler characteristic heuristic for approximating the
  distribution of the largest eigenvalue of an orthogonally invariant random
  matrix}, \bibinfo{journal}{Journal of Statistical Planning and Inference}
  \bibinfo{volume}{138} (\bibinfo{year}{2008}) \bibinfo{pages}{3357--3378}.
\bibitem[{Kuriki and Takemura(2009)}]{kuriki-takemura-2009}
\bibinfo{author}{S.~Kuriki}, \bibinfo{author}{A.~Takemura},
  \bibinfo{title}{volume of tubes and the distribution of the maximum of a
  {G}aussian random field, selected papers on probability and statistics},
  \bibinfo{journal}{American Mathematical Society Translations Series 2}
  \bibinfo{volume}{227} (\bibinfo{year}{2009}) \bibinfo{pages}{25--48}.
\bibitem[{Leykin(2004)}]{Anton2004}
\bibinfo{author}{A.~Leykin}, \bibinfo{title}{Algorithmic proofs of two theorems
  of {S}tafford}, \bibinfo{journal}{Journal of Symbolic Computation}
  \bibinfo{volume}{38} (\bibinfo{year}{2004}) \bibinfo{pages}{1535--1550}.
\bibitem[{Muirhead(2005)}]{Muirhead2005}
\bibinfo{author}{R.~J. Muirhead}, \bibinfo{title}{Aspects of multivariate
  statistical theory}, \bibinfo{publisher}{Wiley}, \bibinfo{year}{2005}.
\bibitem[{Selberg(1944)}]{selberg}
\bibinfo{author}{A.~Selberg}, \bibinfo{title}{Remarks on a multiple integral},
  \bibinfo{journal}{Norsk Matematisk Tidsskrift} \bibinfo{volume}{26}
  (\bibinfo{year}{1944}) \bibinfo{pages}{71--78}.
\bibitem[{Takayama et~al.(2019)Takayama, Jiu, Kuriki, Takayama and
  Zhang}]{ElectronicYZ}
\bibinfo{author}{N.~Takayama}, \bibinfo{author}{L.~Jiu},
  \bibinfo{author}{S.~Kuriki}, \bibinfo{author}{N.~Takayama},
  \bibinfo{author}{Y.~Zhang}, \bibinfo{title}{Supplementary electronic material
  to the article computations of the expected {E}uler characteristic for the
  largest eigenvalue of a real non-central {W}ishart matrix},
  \bibinfo{year}{2019}.
  \bibinfo{note}{\href{https://yzhang1616.github.io/ec1/ec1.html}{https:/$\!$/yzhang1616.github.io/ec1/ec1.html}}.
\bibitem[{Takemura and Kuriki(1999)}]{takemura-kuriki-1999}
\bibinfo{author}{A.~Takemura}, \bibinfo{author}{S.~Kuriki},
  \bibinfo{title}{Shrinkage to smooth non-convex cone: {P}rincipal component
  analysis as {S}tein estimation}, \bibinfo{journal}{Communications in
  Statistics: Theory and Methods} \bibinfo{volume}{28} (\bibinfo{year}{1999})
  \bibinfo{pages}{651--669}.
\bibitem[{Taylor and Worsley(2013)}]{taylor-worsley-2013}
\bibinfo{author}{J.~E. Taylor}, \bibinfo{author}{K.~J. Worsley},
  \bibinfo{title}{Detecting sparse cone alternatives for {G}aussian random
  fields, with an application to f{MRI}}, \bibinfo{journal}{Statistica Sinica}
  \bibinfo{volume}{23} (\bibinfo{year}{2013}) \bibinfo{pages}{1629--1656}.
\bibitem[{Worsley(1995)}]{worsley-1995}
\bibinfo{author}{K.~J. Worsley}, \bibinfo{title}{Boundary corrections for the
  expected {E}uler characteristic of excursion sets of random fields, with an
  application to astrophysics}, \bibinfo{journal}{Advances in Applied
  Probability} \bibinfo{volume}{27} (\bibinfo{year}{1995})
  \bibinfo{pages}{943--959}.
\bibitem[{Zeilberger(1991)}]{Zeilberger1991}
\bibinfo{author}{D.~Zeilberger}, \bibinfo{title}{The method of creative
  telescoping}, \bibinfo{journal}{Journal of Symbolic Computation}
  \bibinfo{volume}{11} (\bibinfo{year}{1991}) \bibinfo{pages}{195--204}.
\bibitem[{Zhang(2015)}]{woods-zhang-2015}
\bibinfo{author}{L.~Zhang}, \bibinfo{title}{volumes of orthogonal groups and
  unitary groups}, \bibinfo{year}{2015}. \bibinfo{note}{ArXiv:1509.00537}.

\end{thebibliography}


\end{document}